\documentclass{article}
\usepackage{amsmath,mathtools}

\newcommand{\xddots}{%
  \raise 4pt \hbox {.}
  \mkern 6mu
  \raise 1pt \hbox {.}
  \mkern 6mu
  \raise -2pt \hbox {.}
}

\usepackage{alltt}
\usepackage[vcentermath]{youngtab}
\usepackage[all]{xy}
\usepackage{latexsym,amssymb,amsmath,amsfonts, amsthm, xcolor}

\newcommand{\RR}{\mathbb{R}}

\def\vbar{\mathchoice{\vrule height6.3ptdepth-.5ptwidth.8pt\kern-.8pt}
  {\vrule height6.3ptdepth-.5ptwidth.8pt\kern-.8pt}
  {\vrule height4.1ptdepth-.35ptwidth.6pt\kern-.6pt}
  {\vrule height3.1ptdepth-.25ptwidth.5pt\kern-.5pt}}
\def\fudge{\mathchoice{}{}{\mkern.5mu}{\mkern.8mu}}
\def\bbc#1#2{{\rm \mkern#2mu\vbar\mkern-#2mu#1}}
\def\bbb#1{{\rm I\mkern-3.5mu #1}}
\def\bba#1#2{{\rm #1\mkern-#2mu\fudge #1}}
\def\bb#1{{\count4=`#1 \advance\count4by-64 \ifcase\count4\or\bba A{11.5}\or
 \bbb B\or\bbc C{5}\or\bbb D\or\bbb E\or\bbb F \or\bbc G{5}\or\bbb H\or
  \bbb I\or\bbc J{3}\or\bbb K\or\bbb L \or\bbb M\or\bbb N\or\bbc O{5} \or
  \bbb P\or\bbc Q{5}\or\bbb R\or\bbc S{4.2}\or\bba T{10.5}\or\bbc U{5}\or
   \bbb P\or\bbc Q{5}\or\bbb R\or\bba S{8}\or\bba T{10.5}\or\bbc U{5}\or
  \bba V{12}\or\bba W{16.5}\or\bba X{11}\or\bba Y{11.7}\or\bba Z{7.5}\fi}}

\newcommand{\vs}{\vspace{0.25cm}}

\newtheorem{theorem}{Theorem}
\newtheorem{itlemma}{Lemma}[section]
\newtheorem{itproposition}[itlemma]{Proposition}
\newtheorem{itcorollary}[itlemma]{Corollary}
\newtheorem{itremark}[itlemma]{Remark}
\newtheorem{itremarks}[itlemma]{Remarks}
\newtheorem{itdefinition}[itlemma]{Definition}
\newtheorem{itexample}[itlemma]{Example}

\newenvironment{lemma}{\begin{itlemma}\rm}{\end{itlemma}} 
\newenvironment{remark}{\begin{itremark}\rm}{\end{itremark}} 
\newenvironment{remarks}{\begin{itremarks} \rm}{\end{itremarks}}
\newenvironment{corollary}{\begin{itcorollary}\rm}{\end{itcorollary}}
\newenvironment{proposition}{\begin{itproposition}\rm}{\end{itproposition}}
\newenvironment{definition}{\begin{itdefinition}\rm}{\end{itdefinition}}
\newenvironment{example}{\begin{itexample}\rm}{\end{itexample}}
\newenvironment{fact}{\noindent {{\bf Fact}}:\ \ }{\hfill \medskip}
\newenvironment{claim}{\noindent {\em Claim}. \ \ }{\hfill \medskip}
\newcommand{\be}[1]{\begin{equation}\label{#1}}
\newcommand{\ee}{\end{equation}}
\newcommand{\bl}[1]{\begin{lemma}\label{#1}}
\newcommand{\br}[1]{\begin{remark}\label{#1}}
\newcommand{\brs}[1]{\begin{remarks}\label{#1}}
\newcommand{\bt}[1]{\begin{theorem}\label{#1}}
\newcommand{\bd}[1]{\begin{definition}\label{#1}}
\newcommand{\bp}[1]{\begin{proposition}\label{#1}}
\newcommand{\bc}[1]{\begin{corollary}\label{#1}}
\newcommand{\bfact}[1]{\begin{fact}\label{#1}}
\newcommand{\bex}[1]{\begin{example}\label{#1}}
\newcommand{\ec}{\end{corollary}}
\newcommand{\efact}{\end{fact}}
\newcommand{\eex}{\end{example}}
\newcommand{\el}{\end{lemma}}
\newcommand{\er}{\end{remark}}
\newcommand{\ers}{\end{remarks}}
\newcommand{\et}{\end{theorem}}
\newcommand{\ed}{\end{definition}}
\newcommand{\ep}{\end{proposition}}
\newcommand{\epr}{\end{proof}}
\newcommand{\bpr}{\begin{proof}}
\newcommand{\bcl}{\begin{claim}}
\newcommand{\ecl}{\end{claim}}

\newcommand{\bi}{\begin{itemize}}
\newcommand{\ei}{\end{itemize}}
\newcommand{\ben}{\begin{enumerate}}
\newcommand{\een}{\end{enumerate}}

\usepackage{graphicx}

\title{Sub-Riemannian Geodesics on $SL(2, \mathbb{R})$ }
\author{Domenico D'Alessandro\thanks{{ Department of Mathematics, Iowa State University, Ames IA 50011, U.S.A., daless@iastate.edu}} \and Gunhee Cho\thanks{Department of Mathematics, University of California Santa Barbara, Santa Barbara CA 93106-3080, U.S.A. gunheecho@ucsb.edu}}
\date{\today}

\begin{document}

\maketitle 
\begin{abstract}
We explicitly describe  the length minimizing geodesics for a sub-Riemannian structure of the elliptic type defined on $SL(2, \mathbb{R})$.  Our method uses a symmetry reduction which translates the problem into a Riemannian problem on a two dimensional quotient space,  on which projections of geodesics can be easily visualized. As a byproduct, we obtain an alternative derivation of the characterization of the cut-locus obtained in \cite{BoscaRossi}. We use  classification results for  three dimensional right invariant sub-Riemannian structures on Lie groups  \cite{AGBD}, \cite{Biggs},  \cite{HB2} to identify exactly   automorphic  structures on which our results apply. 
\end{abstract}

\vs 

\vs
\vs

\vs 

\vs
\vs

\section{Introduction}

A {\bf sub-Riemannian structure} is a triple $(M,\Delta,g)$, where $M$ is a connected smooth manifold, $\Delta$ is a constant rank distribution\footnote{That is,  a smooth assignment to every $x \in M$ of a subspace $\Delta_x \subseteq T_xM$ of the tangent space $T_xM$ at $x$ of constant dimension.} and $g$ a smooth metric defined on $\Delta$.   The distribution $\Delta$ is usually described as the span of a set of vector fields $\{f_1,...,f_m\}$ ({\it a frame}), that is, 
$\Delta=\texttt{span}\{f_1,...,f_m\}$   and it is assumed to be nonintegrable, that is,  $[ \Delta, \Delta ] \nsubseteq \Delta$. The vector fields $\{f_1,...,f_m\}$   are also assumed to be orthonormal with respect to the  (sub-Riemannian) metric $g$, defined on $\Delta$,  i.e., $g_{p}(f_j(p),f_k(p))=\delta_{j,k}$. If $\Delta=TM$ nonintegrability is clearly not possible and one recovers the case of a Riemannian manifold. 
  A {\it horizontal curve} $\gamma: [a,b] \subseteq R \rightarrow M $  is a  Lipschitz continuous curve satisfying, a.e., $\dot \gamma(t) \in \Delta_{\gamma(t)}$ and {\it minimizing geodesics}  are horizontal curves which minimize the {\it sub-Riemannian distance}   between two points $p$ and $q$ in $M$ (also called the {\it Carnot-Caratheodory distance})  where such  a distance is defined by (with $\gamma(0)=p$, $\gamma(T)=q$) 
$$
d(p,q)=\min_{\gamma} \int_0^T\sqrt{g_{\gamma(t)}(\dot \gamma(t), \dot \gamma(t))}dt, 
$$ 
with the minimum taken over all the horizontal curves $\gamma$ joining $p$ and $q$. Such a minimum  exists, under the assumption that the frame $\{f_1,...,f_m\}$ is {\it bracket generating}, i.e., the Lie algebra generated by these vector fields has full dimension.This is the classical {\it Chow-Raschevski theorem}  (see, e.g., \cite{ABB}), and this will  be the case in the problem we shall treat here.  The problem of finding sub-Riemannian geodesics is equivalent (in that it has the same solution) to problems in control theory, where one considers the system 
$$
\dot x=\sum_{j=1}^q u_j(t) f_{j},
$$ 
and wants to find control functions $\{ u_j\}$ in order to drive the state $x$ from  $x(0)=p$ to $x(T)=q$, in minimum time subject to a $L_{\infty}$ bound on the norm of the control or, equivalently, with fixed time $T$, and minimizing an energy type of functional $\int_0^T \sum_j u_j^2(t)dt$ (see, e.g., \cite{AD}).  In particular, sub-Riemannian geodesics $\gamma$  parametrized by arclength ($\| \dot \gamma \|=1$, a.e.)  coincide with minimum time trajectories when the norm of the control is bounded by $1$ \cite{ABB} and, in fact,  the minimum time can be taken as the sub-Riemannian distance. In view of this fact, in the following discussion,  we shall use the words {\it minimum time} and {\it distance}  interchangeably.  Methods of geometric control theory \cite{Agrachev}, such as 
the {\it Pontryaging Maximum Principle}  are used to solve sub-Riemannian problems.

Of special importance for the methods we shall use in this paper is the concept of isometry. A (global) {\bf isometry} between two sub-Riemannian structures $(M,\Delta, g)$ and $(M^{'}, \Delta^{'}, g^{'})$ is  a diffeomorphism $\phi \, : \, M \rightarrow M^{'}$ which preserves both the distributions $\Delta$, $\Delta^{'}$ and the metrics $g$, $g^{'}$, that is,  $\phi_* \Delta=\Delta^{'}$ and $g=\phi^*g^{'}$. 
It follows from the definition that isometries map minimizing geodesics to minimizing geodesics.  An important special case is when the isometry $\phi$ maps a sub-Riemannian structure to itself. In this case, 
the set of all such isometries form a group $Iso(M, \Delta, g)$  Families of geodesics that can be mapped one to the other by elements of such group can be identified as representing a curve in  $M/Iso (M, \Delta, g)$,  the quotient space of $M$ under the action of the transformation group $Iso(M, \Delta, g)$, and, in fact, the problem of finding the geodesics can be treated on the quotient space   $M/Iso(M, \Delta, g)$. This is the essence of the method of {\it symmetry reduction} \cite{AD} which will be used in this paper. For a 
comprehensive introduction to sub-Riemannian geometry we refer to \cite{ABB}.

This paper deals with a sub-Riemannian structure on $SL(2)=SL(2, \mathbb{R})$, the Lie group of area preserving transformations on $R^2$. For this structure, we shall give a complete and explicit description of the {\it optimal synthesis}, i.e., the description of all the geodesics. For Lie groups, the frame of vector fields $\{ f_1,...,f_m\}$ defining the distribution $\Delta$ is usually taken to be made  of right (or left) invariant vector fields, i.e., for every $a \in M$, $R_{a*} f_j (x)=f_{j}(R_a x)$ ($L_{a*} f_j (x)=f_{j}(L_a x)$), for every $x \in M$ and $j=1,...,m$, where $R_a$ ($L_a$) denotes the right (left) translation on the group $M$ by $a$. This way the distribution $\Delta$ is described by a set of elements in the tangent space at the identity ${\bf 1}$,  $T_{\bf 1} M$,  identified with the {\it Lie algebra} of the Lie group. We shall focus here on right invariant structures but similar things can be said for left invariant ones. The metric $g$ is also assumed to be right  invariant, that is, for $X_p, Y_p \in T_{p}M$ and 
$a \in M$, $g_{R_a(p)}( R_{a*} X_p, R_{a*}Y_p)= g_p(X_p, Y_p)$. Therefore the metric can be defined by giving an inner product on the portion of the Lie algebra $T_{\bf 1} M$ corresponding to the distribution $\Delta$ at ${\bf 1}$   and it is then extended  by right invariance at any point $p \in M$. For right invariant structures the right multiplication by $a$, $R_a$,  is an isometry 
$R_{a*} \, : \,  (M, \Delta, g) \rightarrow (M, \Delta, g)$  and,  in particular, it maps minimizing geodesics to minimizing geodesics. Therefore, there is no loss of generality in studying only the geodesics starting from the identity ${\bf 1}$, because, in general,  the geodesic minimizing the length between $p$ and $q$ is obtained as $R_p \circ \gamma$, where $\gamma$ is the length minimizing geodesic from ${\bf 1}$ to $qp^{-1}$. In view of this fact,  in the sequel we shall only consider geodesics with the initial point given by the identity. In applying symmetry reduction,  one considers therefore the subgroup of $Iso(M, \Delta,g)$, of isometries $\phi \, : \, (M,\Delta,g) \rightarrow (M,\Delta,g)$ which leave the identity of the group $M$ unchanged  i.e., $\phi({\bf 1}) ={\bf 1}$. This subgroup will be denoted by $Iso_{\bf 1}(M, \Delta, g)$. If $Iso_{\bf 1} (M, \Delta, g)$ is a continuous group, the quotient space $M/Iso_{\bf 1} (M, \Delta, g)$ will have dimension strictly less than the dimension of $M$. In particular, if $Iso_{\bf 1} (M, \Delta, g)$ is a compact connected Lie group, it follows from the theory of Lie transformation groups (see, e.g., \cite{Bredon}) that $M/Iso_{\bf 1} (M, \Delta, g)$ is a {\it stratified space}  with one particular  stratum, corresponding to points with minimal  isotropy type,  which is a connected (open and dense) manifold. This is called the {\it regular part} of   $M/Iso_{\bf 1} (M, \Delta, g)$ and, if the minimal isotropy group is discrete,   it has dimension equal to $\dim(M) - \dim\left( Iso_{\bf 1} (M, \Delta, g) \right)$. The remaining part of   $M/Iso_{\bf 1} (M, \Delta, g)$  is called the {\it singular part}.  In applying symmetry reduction,  the regular part of   the quotient space $M/Iso_{\bf 1} (M, \Delta, g)$ may be given a Riemannnian metric so that Riemannian geodesics in $M/Iso_{\bf 1} (M, \Delta, g)$ are projections of sub-Riemannian geodesics in $M$  with the same length \cite{DS} \cite{S}. This is the approach we will follow in this paper to describe {\it all}  the sub-Riemannian geodesics, i.e., the {\it complete optimal synthesis} on $SL(2)$. 

In physics, $SL(2)$ is referred to as {\it 2+1 dimensional anti-de Sitter space}, since it  is  regarded as a three dimensional hypersurface sitting inside a four dimensional space, with  a natural Lorentzian metric of  signature $(2, 1)$. In general, it is a difficult problem to connect  two arbitrary points by Lorentzian (semi- Riemannian) geodesics, (see, e.g., \cite{DCCIMAV}). In sub-Riemannian geometry, the anti-de Sitter fibrations and the Hopf fibrations viewed as sub-Riemannian manifolds have been studied steadily. In particular, by applying the free action from the totally geodesic fiber of these two fibrations, the radial part of the sub-Laplacian can be calculated explicitly (e.g., \cite{FBGC1,FBGC2}). 

The paper is organized as follows. We start by describing the rank 2 sub-Riemannian structure $(M, \Delta, g)$ on $SL(2)$ we want to study (section \ref{SRSSL2}) and in section \ref{SR9} we pick $SO(2)$   as a Lie subgroup  of $Iso_{\bf 1}(M, \Delta, g)$ which we plan to use to carry out  symmetry reduction.\footnote{The full $Iso_{\bf 1}$ group in this case is $O(2)$ according to the results of \cite{Biggs} but we will take only the connected component containing the identity in order to keep in line with the  standard theory of Lie transformation groups \cite{Bredon}. This results in an extra symmetry in the plot of the optimal synthesis in figure \ref{Figura1} (symmetry about the $x$ axis) which however does not complicates or invalidate our analysis.} The quotient space $SL(2)/SO(2)$ is described in section \ref{Desk}. It is given by the part of the $R^2$ plane outside a disk of radius one and including  the circle of radius one. The regular part corresponds to the part strictly outside the disk. Our goal is to arrive at a complete description of all the sub-Riemannian geodesics in terms of certain Riemannian geodesics in the regular part of $SL(2)/SO(2)$. These are the images under the natural projection $\pi \, : \,  SL(2) \rightarrow SL(2)/SO(2)$ of sub-Riemannian geodesics. In order to achieve this,  we have to define a Riemannian metric on the regular part of $SL(2)/SO(2)$ so that the natural projection $\pi$ preserves the lengths. In section \ref{VFQS} we describe the tangent  vector in $T SL(2)/SO(2)$ obtained by projecting the tangent vectors defining the sub-Riemannian structure on $SL(2)$ at any point. By imposing that these tangent vectors  are orthonormal we obtain the appropriate Riemannian metric in section \ref{RM}, which leads to writing down the geodesic equations in section \ref{GEQS}. The geodesics in the quotient space are obtained in section \ref{GQS}  by projecting the sub-Riemannian geodesics which were obtained in \cite{BoscaRossi}. These curves in fact satisfy the geodesic equations as verified in Appendix A. 
Visualizing the geodesics in the (two dimensional) plane allows us to easily select the geodesic joining the identity to any desired final point. This is the main contribution of this work which is described in sections \ref{GDOS} and \ref{Exe}. In these sections,  we give the {\it complete optimal synthesis}, that is, the description of 
all the geodesics and show how to solve the length minimization problem between any two points in practice. As a 
byproduct,   we obtain an alternative proof of the characterization of the cut-locus which was first obtained in \cite{BoscaRossi}. The results of \cite{BoscaRossi} are  the most related to the ones presented in this paper. As compared to their characterization of the geodesics, the description given here is more explicit and constructive, giving a practical  way to solve the length minimizing problem between any two points.   Since the symmetry reduction technique was also used  for the solution of a sub-Riemannian problem on $SU(2)$ (\cite{Adpre}), and  the Lie algebras $sl(2)$ and $su(2)$ are both real forms of the same complex Lie algebra $sl(2,C)$, it is natural to ask about the relation between the two problems. We do that in section \ref{Relat5}. On any Lie group $M$ we can define many sub-Riemannian structures even if we require right invariance. Structures $(M,\Delta, g)$ and $(M, \Delta^{'}, g^{'})$ which are related by an automorphism $\phi$ of the Lie group $M$, i.e., $\phi_* \Delta=\Delta^{'}$ $g=\phi^* g^{'}$, are seen as equivalent, and sub-Riemannian geodesics  in one structure are mapped  to sub-Riemannian geodesics in the other structure. With the help of the classification of right invariant sub-Riemannian structures on three dimensional Lie groups \cite{AGBD}, \cite{Biggs}, \cite{HB2}, we identify in section  \ref{Ident5} the structures on $SL(2)$ which are equivalent via automorphism to the one we have treated which leaves as an open problem to characterize the optimal synthesis on the remaining sub-Riemannian structures on the Lie group $SL(2)$.

\section{Sub-Riemannian structure on $SL(2)$}\label{SRSSL2}
 On $SL(2)$ we consider the three right invariant vector fields corresponding to the elements of the Lie algebra $sl(2)$, $\{A_0, A_1,A_2\}$,  defined by 
\be{f0f1f2}
A_0:=\frac{1}{2}\begin{pmatrix} 0 & -1 \cr 1 & 0  \end{pmatrix} , \qquad 
A_1:=\frac{1}{2}\begin{pmatrix} 0 & 1 \cr 1 & 0  \end{pmatrix} , \qquad 
A_2:=\frac{1}{2}\begin{pmatrix} 1 & 0 \cr 0 & -1  \end{pmatrix} , 
\ee
and satisfying the commutation relations 
\be{commurel}
[A_0, A_1]=-A_2, \qquad [A_0, A_2]=A_1, \qquad [A_1,A_2]=A_0. 
\ee
The basis $\{A_0, A_1,A_2\}$ gives a {\it Cartan decomposition} (see, e.g.,  \cite{Helgason}) of $sl(2)$ in that by setting ${\cal K}=\texttt{span} \{ A_0 \}$, ${\cal P}=\texttt{span} \{ A_1, A_2\}$, we have from (\ref{commurel}) 
\be{Cartancomm}
[{\cal K}, {\cal K} ]\subseteq {\cal K}, \qquad [{\cal K}, {\cal P}] \subseteq {\cal P}, \qquad [{\cal P}, {\cal P}] \subseteq {\cal K}. 
\ee
We shall denote by, $f_0,$ $f_1,$ $f_2$, the right invariant vector fields corresponding respectively to $A_0,$ $A_1,$ and $A_2$, and by $f_{0,1,2}(X)$ their values at 
$X \in SL(2)$.  
These vector fields  are orthonormal at every point $X \in SL(2)$  with the right invariant  metric given by (the definition on the Lie algebra $sl(2)$)\footnote{Up to a proportionality factor,  this inner product coincides with the {\it Killing metric} defined on the Lie algebra $sl(2)$ by  
$Kill(B,C)=Tr(ad_B ad_C)$ where $ad_A$ is the adjoint representation of $A$  which calculated on $\texttt{span}\{ A_0, A_1, A_2\}$ (cf. (\ref{addav})) which  is negative definite on $\texttt{span}\{ A_0\}$ and positive definite  on $\texttt{span} \{ A_1, A_2\}$.} 
\be{metrica}
g(B,C):=2Tr(BC^T). 
\ee 
The sub-Riemannian structure $\Delta$ on $SL(2)$ we shall study is given by the sub-bundle of the tangent bundle  determined by $\{A_1,A_2\}$ with the sub-Riemannian metric given by the restriction of  the metric $g$ in (\ref{metrica}) to $\Delta(X) \subset T_XSL(2)$. We want to describe sub-Riemannian minimizing geodesics connecting any two points in $SL(2)$. Because of the right invariance of the distribution, it is enough to study geodesics starting from the identity.    

\section{Symmetry Reduction}\label{SR9} 
The group $SO(2)$ acting on $SL(2)$ by conjugation ($X \in SL(2)$, $K \in SO(2)$,  $X \rightarrow KXK^T$) is a group of {\it isometries} for the above sub-Riemannian structure which leaves the identity fixed,  that is,  $SO(2) \subseteq Iso_{\bf 1} (SL(2), \Delta, g)$, since $K\left( \texttt{span}\{ A_1, A_2\} \right) K^T = \texttt{span}\left\{ A_1, A_2 \right \} $. In particular,  it preserves the sub-Riemannian length of curves starting from the identity and  if $\gamma=\gamma(t)$ is a sub-Riemannian geodesic connecting ${\bf 1}$ to $X_f \in SL(2)$, $K\gamma K^T=K\gamma(t) K^T$ is a sub-Riemannian geodesic connecting 
${\bf 1}$ to $KX_fK^T \in SL(2)$, with the same sub-Riemannian length. The idea of {\it{symmetry reduction}} is  therefore  to consider the quotient space under such an action, $SL(2)/SO(2)$,  with a Riemannian metric such that the Riemannian distance between two points in  $SL(2)/SO(2)$ is equal  to the sub-Riemannian distance between  two corresponding points in $SL(2)$. In particular, minimizing Riemannian geodesics in $SL(2)/SO(2)$ will correspond to minimizing sub-Riemannian geodesics in $SL(2)$ (\cite{AD}, \cite{DS}, \cite{S}). We denote by $\pi$ the natural projection 
\be{pi}
\pi \, : \, SL(2) \rightarrow SL(2)/SO(2). 
\ee
In the next section we describe the quotient space $SL(2)/SO(2)$. The space  $SL(2)/SO(2)$ is not a manifold but it is a stratified space which contains a subset (a stratum) which is connected and open and dense in 
$SL(2)/SO(2)$, called  the {\it regular part}   of the quotient  $SL(2)/SO(2)$, and which can be given the structure of a manifold. It is on this regular part that we will define a Riemannian metric. According to standard results in Lie transformation groups \cite{Bredon}, the regular part is the image under the natural projection $\pi$ of the set of points in $SL(2)$ which has minimal isotropy group.\footnote{In general, for a Lie transformation group $G$ on a manifold $M$, one considers  an equivalence relation between subgroups 
$H_1$ and $H_2$ of $G$ saying that $H_1$ is equivalent to $H_2$ if there exists a $k \in G$ such that 
$kH_1 k^{-1}=H_2$.  Denote the equivalence class containing $H$ by $(H)$. Some of these equivalence classes contain subgroups which are isotropy groups of some elements of $M$. These are called {\it isotropy types}. On isotropy types, one can consider a partial ordering by saying that $(H_1) \leq (H_2)$ if $H_1$ is conjugate (via an element of $G$)  to a subgroup of $H_2$. According to a basic result in the theory of Lie transformation groups, there exists a minimal isotropy type in this ordering which corresponds to a stratum in $M/G$ which is a connected, open and dense manifold in $M/G$.  This is called the {\it regular part} of the quotient space $M/G$ the remaining part is the {\it singular part}. A more detailed discussion can be found in \cite{AD} and a complete treatment in \cite{Bredon}.}

\section{Description of $SL(2)/SO(2)$}\label{Desk} 
The equivalence relation $\sim$  defining $SL(2)/SO(2)$ is such that, given   
$$
X_1:=\begin{pmatrix} a_1 & b_1 \cr c_1 & d_1  \end{pmatrix}, \qquad 
X_2:=\begin{pmatrix} a_2 & b_2 \cr c_2 & d_2  \end{pmatrix},
$$
$X_1 \sim X_2$ if and only if $a_1+d_1=a_2+d_2$ and $b_1-c_1=b_2-c_2$. To see this notice that for a matrix 
\be{genmatX}
X=\begin{pmatrix} a & b \cr c & d  \end{pmatrix}, 
\ee
the two quantities 
\be{xey}
x:=\frac{a+d}{2},  \qquad y:=\frac{b-c}{2}
\ee
are invariant under the action of $K \in SO(2)$, $x$ because it is $\frac{1}{2}$ the trace of $X$,  and $y$ because, we can write 
$$
X=\begin{pmatrix} a & \frac{b+c}{2} \cr 
\frac{b+c}{2} & d \end{pmatrix}+ \begin{pmatrix} 0 & {y} \cr 
- {y} & 0  \end{pmatrix}, 
$$ 
and 
$$
KXK^T=K\begin{pmatrix} a & \frac{b+c}{2} \cr 
\frac{b+c}{2} & d \end{pmatrix} K^T+ \begin{pmatrix} 0 & {y} \cr 
- {y} & 0  \end{pmatrix}:=  \begin{pmatrix} \tilde a & \tilde k +{y}  \cr 
\tilde k-{y}  & \tilde d  \end{pmatrix}, 
$$
where the invariance is due to $\frac{1}{2}\left( 
\left(\tilde k+{y}\right)-
\left(\tilde k-{y}\right) \right)=y$. Viceversa assume that two matrices $X_1$ and $X_2$ in $SL(2)$ have the same values of $x$ and $y$. Then we have to show that there exists $K\in SL(2)$ such that $X_2=KX_1K^T$. Write 
\be{X1X2}
X_{1,2}=\begin{pmatrix} x & y \cr -y &x \end{pmatrix}+ \begin{pmatrix} m_{1,2} & k_{1,2} \cr  k_{1,2} & -m_{1,2} \end{pmatrix}, 
\ee
and notice that,  from the condition $\det(X_1)=\det(X_2)=1$,  we obtain 
\be{POL}
m_1^2+ k_1^2=m_2^2+k_2^2={x^2} + {y^2}-1.  
\ee
Thus the points $(m_1,k_1)$ and $(m_2,k_2)$ are on the same circle with radius 
${x^2} + {y^2}-1$. Now choose $\theta$ such  that
\be{transfor}
\begin{pmatrix} 
m_2 \cr k_2
\end{pmatrix} =\begin{pmatrix} \cos(2 \theta) & \sin(2 \theta) \cr - \sin(2\theta) & \cos(2\theta) \end{pmatrix} \begin{pmatrix} 
m_1 \cr k_1
\end{pmatrix} .
\ee 
With this choice, a direct calculation shows that using 
\be{kappa}
K=\begin{pmatrix} \cos(\theta) & \sin(\theta) \cr -\sin(\theta) & \cos(\theta) \end{pmatrix}, 
\ee
we have 
$$
KX_1 K^T=X_2. 
$$

The above argument and equation (\ref{POL}) show that the points $(x,y)$ of the $x-y$ plane such that $x^2+y^2 \geq 1$ are in one to one correspondence with  the points of the quotient space 
$SL(2)/SO(2)$ which can  be identified therefore with the portion of the plane {\it outside} the open disc of radius $1$, together with the boundary of the unit disc itself.  Points on the boundary of the disc correspond to points in $SL(2)$ whose isotropy group (see, e.g. \cite{Bredon})  is the full $SO(2)$ group. This is the {\it singular part} of the quotient space, and the corresponding (singular) part on $SL(2)$ correspond to matrices such that $x^2+y^2=1$, that is for (\ref{genmatX}), using $x^2+y^2=\frac{(a+d)^2}{4} +\frac{(b-c)^2}{4}=1$, together with $ad-bc=1$, 
we have $(a-d)^2+(b-c)^2=0$, so that, the matrices in the singular part have the form 
\be{singpart7}
X=\begin{pmatrix} x & y \cr -y & x \end{pmatrix}, 
\ee
with $x^2+y^2=1$
The rest of the plane outside the disk is the {\it regular (open, dense and connected) part} of the quotient space, which we denote by $\left( SL(2)/SO(2) \right)_{reg}$.  If we use the representation (\ref{X1X2}) of a matrix $X\in SL(2)$, $X=\begin{pmatrix} x & y  \cr -y & x \end{pmatrix}+ \begin{pmatrix} m & k \cr k & -m \end{pmatrix}$, matrices corresponding to points in this region have $m^2+ k^2=x^2+y^2-1>0$. If $K\in SO(2)$ in (\ref{kappa}) is in the isotropy group of this element, then using $KXK^T=X$ and specializing (\ref{transfor}), we obtain 
$$
\begin{pmatrix} m \cr k  \end{pmatrix}=\begin{pmatrix} \cos(2 \theta) & \sin(2\theta) \cr -\sin(2\theta) & \cos(2 \theta) \end{pmatrix} \begin{pmatrix} m \cr k  \end{pmatrix}. 
$$ 
Taking the inner product with $\begin{pmatrix} m \cr k  \end{pmatrix}$ gives $\cos(2 \theta)(m^2+k^2)=(m^2+k^2)$, that is,  $\theta=k\pi$. Thus the isotropy group of this type of elements  is made of only $\pm {\bf 1}$. This is the (discrete)  minimal isotropy type.  This  regular part  of the orbit space will be given  the structure of a Riemannian manifold in the next section. The orbit at each point $(x,y)$  is a circle with square radius ${x^2+y^2}-1$.\footnote{We notice, in particular,  the similarity with the $SU(2)$ case  which was treated in \cite{Adpre} where the quotient space under the action of $SO(2)$  is the {\it interior} of the unit disk. The relation with the $SU(2)$ case is discussed in section \ref{Relat5}. }

The $x$ axis ($y=0$) represents the matrices $X$ in $SL(2)$ which are in the same class as their inverse. These are the matrices in $SL(2)$ of the form $X=\begin{pmatrix} a & b \cr b &d \end{pmatrix}$. Their inverse 
$X^{-1}=\begin{pmatrix} d & -b \cr - b &a \end{pmatrix}$ is in the same class  because $\begin{pmatrix} 0 & 1 \cr - 1 &0  \end{pmatrix} X \begin{pmatrix} 0 & -1 \cr  1 &0  \end{pmatrix}=X^{-1}$. Conversely, if $X$ and $X^{-1}$ are in the same class, they must have the same $y$. However the $y$ of $X$ and $X^{-1}$ are opposite to each other (while the $x$'s are the same) which implies $y=0$ and therefore the $(1,2)$ and the $(2,1)$ element of the matrix are the same,  as in  
$X=\begin{pmatrix} a & b \cr b &d \end{pmatrix}$.

\section{Projections of tangent vectors}\label{VFQS}
Consider a point $X$ belonging to $SL(2)$ and the values at this point of the vector fields  defining the given distribution $\Delta$, i.e., $f_1(X)$ and $f_2(X)$. We now calculate,  for every point $\pi(X)$ in the regular part of $SL(2)/SO(2)$, 
$\pi_* f_1(X)$ and $\pi_* f_2(X)$.  We do this in the basis corresponding to the coordinates $x$ and $y$ defined in the previous section, i.e., $\left \{ \frac{\partial}{\partial x}, \frac{\partial}{\partial y} \right \}$.  
We have 
$$
\pi_* f_1(X)=\left( \pi_* f_1(X) x \right) \frac{\partial}{\partial x} +\left( \pi_* f_1(X) y \right) \frac{\partial}{\partial y} 
$$

Using (\ref{genmatX}), (\ref{xey}),  and (\ref{f0f1f2}) we get 
$$
f_1(X)=\frac{1}{2}\left( c \frac{\partial}{\partial a} +
d \frac{\partial}{\partial b}+a \frac{\partial}{\partial c}+ b \frac{\partial}{\partial d} \right).  
$$
Therefore 
$$
\pi_* f_1(X)x=f_1(X)(x \circ \pi)=f_1(X)\left( \frac{a+d}{2}\right)= \frac{1}{4}(c+b), 
$$
and analogously, 
$$
\pi_* f_1(X)y=f_1(X)(y \circ \pi)=f_1(X)\left( \frac{b-c}{2}\right)= \frac{1}{4}(d-a).  
$$
Therefore we have 
\be{pistar1}
\pi_* f_1(X)=\frac{1}{4} \left( (b+c) \frac{\partial}{\partial x} + (d-a) \frac{\partial}{\partial y}  \right). 
\ee
An analogous calculation leads to 
\be{pistar2}
\pi_* f_2(X)=\frac{1}{4} \left( (a-d) \frac{\partial}{\partial x} + (b+c) \frac{\partial}{\partial y}  \right). 
\ee
The expressions  (\ref{pistar1}) and (\ref{pistar2}) show that when $x^2+y^2 > 1$, that is,  in the regular part $\pi_*$,  is an isomorphism 
$\pi_*  \, : \, \Delta(X)  \rightarrow T_{\pi(X)} SL(2)/SO(2)$.\footnote{The map is non singular iff $(a-d)^2+ (b+c)^2>0$, which is true if and only if $a^2+d^2+b^2+c^2 > 0$ which is true (adding $2(ad-bc)=2$ on both sides)  if and only if $(a+d)^2+(b-c)^4 > 4$, which is true, recalling the definitions of $x$ and $y$ (\ref{xey}). } In the following section we define a Riemannian metric on the regular part of $SL(2)/SO(2)$ to make $\pi_*$ an isometry between the sub-Riemannian manifold $SL(2)$ and the Riemannian manifold $SL(2)/SO(2)$.

\section{Riemannian metric on $SL(2)/SO(2)$}\label{RM}
Consider the (diagonal) metric $g_Q$  on the regular part of the quotient space  given in the $(x,y)$ coordinates by (at the point $(x,y)$) , 
\be{defimetric1}
g_Q\left( \frac{\partial}{\partial x},  \frac{\partial}{\partial x} \right)=g_Q\left( \frac{\partial}{\partial y},  \frac{\partial}{\partial y} \right)=\frac{4}{x^2+y^2-1}, 
\ee
\be{defimetric2}
g_Q\left( \frac{\partial}{\partial x},  \frac{\partial}{\partial y} \right)=g_Q\left( \frac{\partial}{\partial y},  \frac{\partial}{\partial x} \right)=0.  
\ee
A direct verification using (\ref{pistar1}) and (\ref{pistar2}) shows that 
\be{exem12}
g_Q(\pi_* f_j, \pi_* f_k)=\delta_{j,k}.  
\ee
This means that, since $f_1$ and $f_2$ are orthonormal  in the original sub-Riemannian metric (\ref{metrica}),   $\pi_*$ is an isometry from the sub-Riemannian  manifold $(SL(2), \Delta, g)$ to the Riemannian manifold $\left( SL(2)/SO(2) \right)_{reg}$ with the  metric $g_Q$. As an  example of verifying (\ref{exem12}), we have, using $ad-bc=1$, 
$$
g_Q(\pi_* f_2, \pi_* f_2)=\frac{1}{4(x^2+y^2-1)} \left[ (a-d)^2+(b+c)^2 \right]=
\frac{1}{4(x^2+y^2-1)} \left[ a^2+d^2+b^2+c^2-2 +2(ad-bc)-2 \right]=1,
$$
using (\ref{xey}).

\section{Geodesic equations on $SL(2)/SO(2)$} \label{GEQS}
To compute the geodesic equations associated with the given metric (\ref{defimetric1}) (\ref{defimetric2}) in $\left( SL(2)/SO(2) \right)_{reg}$,  
we first compute the Christoffel symbols for the Levi-Civita connection associated to the given metric with the classical formulas given for example in \cite{DoCarmo} (formula (10) Chapter 2). The result is 
$$
\frac{x}{x^2+y^2-1}=-\Gamma_{x,x}^x= \Gamma_{y,y}^x=-\Gamma_{x,y}^y=-\Gamma_{y,x}^y, 
$$
$$
\frac{y}{x^2+y^2-1}=-\Gamma_{x,y}^x=-\Gamma_{y,x}^x=\Gamma_{x,x}^y=-\Gamma_{y,y}^y. 
$$
The geodesic equations are given in formula (1) of Chapter 3 of \cite{DoCarmo}. They are in our case: 
\be{geodeq1}
\frac{d^2 x}{dt^2}= \frac{x}{x^2+y^2-1} \left( \left(\frac{dx}{dt} \right)^2-\left(\frac{dy}{dt} \right)^2 \right)+\left(\frac{2y}{x^2+y^2-1} \right) \frac{dx}{dt} \frac{dy}{dt} 
\ee
$$
\frac{d^2 y}{dt^2}= \frac{y}{x^2+y^2-1} \left( \left(\frac{dy}{dt} \right)^2-\left(\frac{dx}{dt} \right)^2 \right)+\left(\frac{2x}{x^2+y^2-1}\right) \frac{dx}{dt} \frac{dy}{dt} . 
$$

\section{Geodesics in $SL(2)/SO(2)$}\label{GQS}

The sub-Riemannian geodesics on $SL(2)$ are such that, when projected onto  the quotient space $SL(2)/SO(2)$,  (locally) minimize the distance with respect to the metric described in the previous section. Therefore the projections of the sub-Riemannian geodesics have to be geodesics for the metric we have defined in section \ref{RM}. They have to satisfy equations (\ref{geodeq1}).

The sub-Riemannian geodesics in $SL(2)$ were  calculated in \cite{BoscaRossi} (subsection 3.3.1) using general results on sub-Riemannian manifold having a $K-P$ structure, that is, whose defining distribution corresponds to the ${\cal P}$ part of a Cartan decomposition as defined in (\ref{Cartancomm}).    The sub-Riemannian geodesics  (parametrized by arclength) are of the form (cf. (\ref{f0f1f2}))
\be{SRgeod}
X(t)=e^{(cA_0 +P)t}e^{-cA_0 t}, 
\ee
where $P$ is a linear combination of $A_1$ and $A_2$ with unit norm, and $c$  a real parameter.  Since we are interested in the projection of these trajectories onto the quotient space,  we can, without loss of generality,  consider, instead of $X(t)$,  $KX(t)K^T$ for an arbitrary $K \in SO(2)$. In particular we can take $P=A_1$. Furthermore assume $X_f:=e^{(cA_0 +P)t}e^{-cA_0 t}$. Then 
\be{inversa}
X_f^{-1}=e^{cA_0 t} e^{-(cA_0+P)t}=e^{(-cA_0+e^{cA_0t} P e^{-cA_0t})t} e^{cA_0 t}.  
\ee
From $\pi(X_f)=(x,y) \leftrightarrow \pi(X_f)=(x,-y)$ and (\ref{inversa}) we see that the projection of the sub-Riemannian geodesic  (\ref{SRgeod}) $\pi(X(t))$ is equal  to the one for $X^{-1}(t)$, 
$\pi\left (X^{-1}(t) \right)$ reflected about the $x$ axis and $\pi(X)$ and $\pi(X^{-1})$ correspond 
to values of $c$ which are opposite of each other. Therefore we can restrict ourselves to $c \geq 0$ in (\ref{SRgeod}) and consider only projections of sub-Riemannian geodesics in the upper half plane. Using 
(\ref{inversa}) and the explicit expression of (\ref{SRgeod}) calculated in \cite{BoscaRossi}, we obtain $\pi(X(t))=(x(t), y(t))$ with, setting $s:=\frac{t}{2}$,\footnote{To simplify the notations, our definitions slightly 
differ from the ones in \cite{BoscaRossi}. In particular our $k_2$ corresponds to their $ck_2$. } 
\be{xt}
x=x(t)=k_1(s) \cos(c s)+k_2(s)\sin(cs), 
\ee 
\be{yt}
y=y(t)=k_1(s)\sin(cs)-k_2(s)\cos(cs), 
\ee
where for $c\geq 0$, the functions $k_1=k_1(s)$ and $k_2=k_2(s)$ are defined as: 
\begin{enumerate}
\item[] For $ c\neq 1$, 
\be{traj1}
\begin{cases}
k_1=\frac{e^{\sqrt{1-c^2}s} +e^{-\sqrt{1-c^2}s}}{2}   \\
k_2=c \frac{e^{\sqrt{1-c^2}s} -e^{-\sqrt{1-c^2}s}}{2\sqrt{1-c^2}}. 
\end{cases}
\ee
In particular, 

\begin{enumerate}

 \item[] For $|c|<1$, 
\be{cass1}
 \begin{cases}
k_1=\cosh(\sqrt{1-c^2} s)  \\
k_2=\frac{c}{\sqrt{1-c^2}}\sinh(\sqrt{1-c^2} s) 
\end{cases}
\ee 
  
 \item[] For $|c|>1$ 
\be{cass2}
\begin{cases}
k_1= \cos(\sqrt{c^2-1} s)  \\
k_2=\frac{c}{\sqrt{c^2-1}}\sin(\sqrt{c^2-1} s) 
\end{cases}
\ee

\end{enumerate}

\item[] For $|c|=1$
\be{traj2}
\begin{cases}
k_1= 1 \\
k_2=\texttt{sign}(c)s
\end{cases}. 
\ee

\end{enumerate}

Notice that for $c=0$ the trajectory is along the $x$ axis, in agreement with the above recalled symmetry.  

A direct but tedious calculation verifies  that the trajectories in (\ref{xt})-(\ref{traj2}) satisfy the geodesic equations (\ref{geodeq1}). For completeness, we present this calculation in  Appendix A.

\section{Geometric description of the optimal synthesis}\label{GDOS}

The word {\it optimal synthesis} refers to a qualitative and graphical description of the 
optimal geodesics (from the identity) to any point. Because of  symmetry reduction  the projections of such geodesics can be plotted in the quotient space, $SL(2)/SO(2)$.  Such geodesics are plotted from the initial point $(1,0)$, which corresponds to the identity in $SL(2)$,   to the point where they  cease to be optimal. The locus where sub-Riemannian geodesics  cease to be optimal is the inverse image under the natural projection of the corresponding locus  in the quotient space \cite{AD}. We shall call this locus the {\bf critical locus}.\footnote{This name is used sometime to describe  the locus where the exponential map is singular.}   Another locus of importance  is the {\bf cut locus}, the locus of points that are reached simultaneously by two 
sub-Riemannian geodesics. Such a locus for $SL(2)$ was described in \cite{BoscaRossi}. In the following, 
 we denote by $R-{\cal CR}$ ($SR-{\cal CR}$)  the `Riemannian'  (`sub-Riemannian') critical locus where the Riemannian geodesics (\ref{xt}), (\ref{yt}) (the sub-Riemannian geodesics (\ref{SRgeod}))    lose their optimality. We have, with the natural projection $\pi$ (cf. \cite{AD}) 
 \be{critpoi}
SR-{\cal CR}=\pi^{-1}(R-{\cal CR}). 
 \ee
As a direct application of a result   proved in \cite{AD} (cf. Corollary 3.6 in that paper), we have  the following: 
 \bp{touch}
 A sub-Riemannian geodesic on the manifold $SL(2)$ which crosses the regular part touches the singular part in a point $p$ and then returns to the regular part loses optimality at $p$. 
 \ep 
 
\begin{figure}[h]
\centerline{\includegraphics[width=4.8in,height=3in]{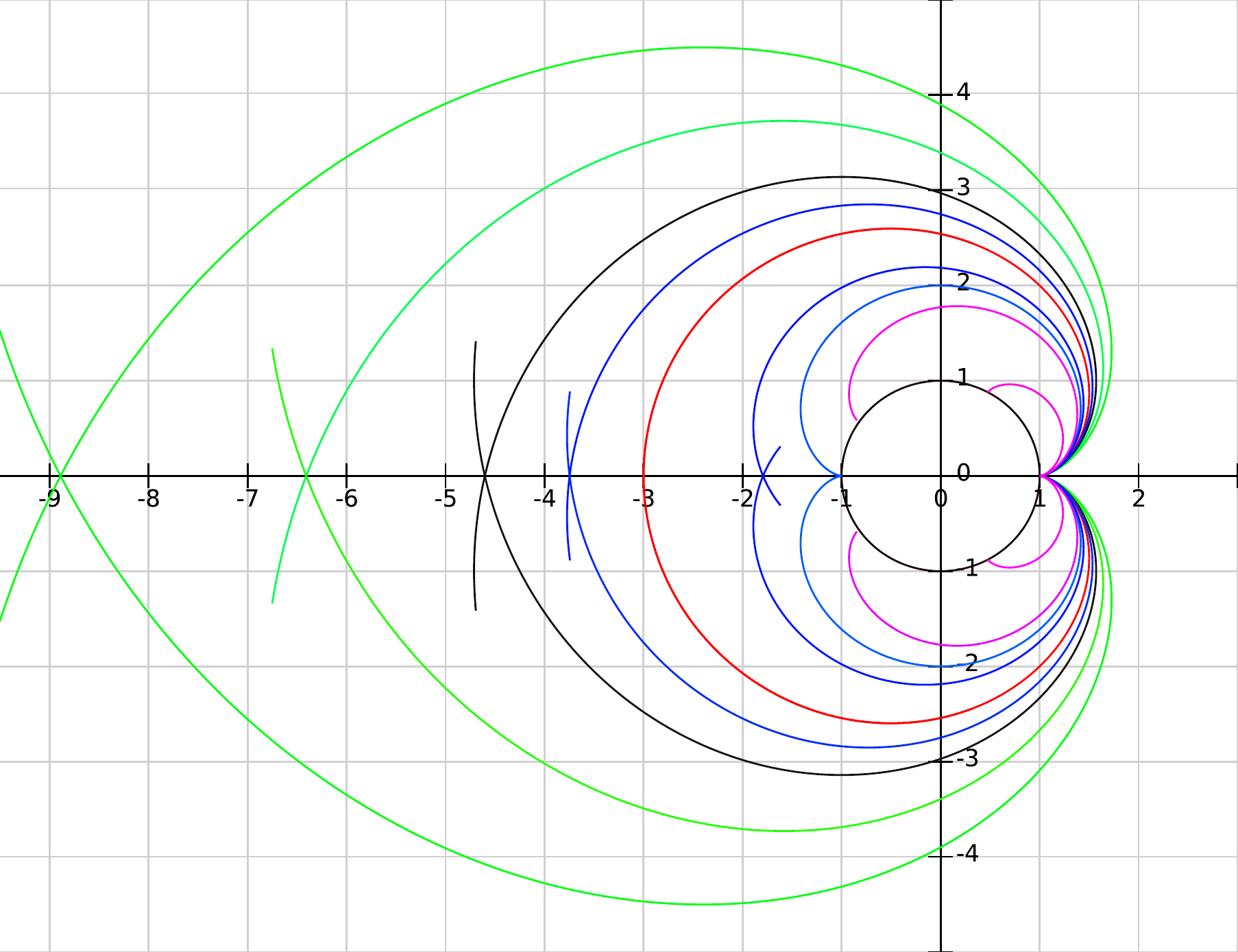}}
\caption{Riemannian geodesics on $SL(2)/SO(2)$ which are projections of sub-Riemannian geodesics in $SL(2)$. All geodesics start from the point $(1,0)$ which is the projection of the identity. The geodesic corresponding to 
$c=0$ follows the $x$ axis in the positive direction. Geodesics corresponding to $|c|<1$ are in green. They intersect the $x$ axis at a point whose $|x|$ value increases and goes to infinity as $|c|\rightarrow 0$. Geodesics corresponding to $|c|=1$ are in black. The geodesics corresponding to $1 < |c| \leq \frac{2}{\sqrt{3}}$ are in blue, except for the ones corresponding to $|c|=\frac{3}{2 \sqrt{2}}$ which are in red. The geodesics corresponding to $c=\pm \frac{2}{\sqrt{3}}$ are the ones that simultaneously reach the point $(-1,0)$. The geodesics corresponding to $|c|> \frac{2}{\sqrt{3}}$ are in purple. They are the ones that `land' on the unit circle, the singular part of the quotient space.}
\label{Figura1}
\end{figure}

We shall use this and the reduction to the quotient space to give a complete description of the optimal synthesis in this section. This is summarized  in Figure \ref{Figura1}.  We shall discuss in the next section how to find the geodesics in practice. for any desired final condition. We remark that symmetry reduction is crucial in this treatment as it allows  easy  visualization of the geodesics in a plane.

As anticipated, the geodesics plotted in Figure \ref{Figura1} are symmetric with respect to the $x$-axis, that is, if $(x_c, y_c)=(x_c(s), y_c(s))$ is the geodesic corresponding to $c$, then 
$(x_{-c}, y_{-c})=(x_c(s),- y_c(s))$. This follows from the fact that $k_1$ ($k_2$) is an even (odd) function of $c$.\footnote{Recall that, as discussed above,  $\pi(X)=(x,y) \leftrightarrow  \pi(X^{-1})=(x,-y)$, that is, reflection about the $x$-axis in the quotient space corresponds to matrix inversion in $SL(2)$.}  Accordingly the geodesics corresponding to  $c=0$, starts from the point $(1,0)$ and moves rightwards along the $x$-axis. 
 We have also that, for every $c$, 
\be{rc2gen}
r_c^2(s):=x_c^2(s)+y_c^2(s)=k_1^2(s)+k_2^2(s)
\ee

\subsection{Geodesics for $|c|<1$}

For $|c|<1$, formula (\ref{rc2gen}) gives 
\be{rc2}
r_c^2(s)=1+\frac{\sinh^2( \sqrt{1-c^2} s)}{1-c^2}.  
\ee
For any $c$, this is an increasing function of $s$. So the geodesics are spirals with increasing distance from the origin. In particular they never self intersect. The geodesic corresponding to $c$ and the one corresponding to $-c$ intersect  when they cross the $x$-axis, with the same value of time $s$. In particular,  notice that if 
$y_c(s)=k_{1,c}(s)\sin(cs)-k_{2,c}\cos(cs)=0$, that is,  $k_{1,c}(s)\sin(cs)=k_{2,c}(s) \cos(cs)$, then $k_{1,-c}(s)\sin(-cs)=-k_{1,c}(s)\sin(cs)=-k_{2,c}(s) \cos(cs)=
k_{2,-c}(s) \cos(-cs)$. Therefore $y_{-c}(s)=0$ and,  by the above recalled symmetry about the $x$ axis,  
$x_{-c}(s)=x_{c}(s)$.

Furthermore the intersections at the $x$-axis of the geodesics 
corresponding to $c$ and $-c$ are the only ones where {\it optimal} geodesics 
with different $c$'s intersect. To see this assume that the geodesic corresponding to $c_1$ and $c_2$ intersect. Because of (time) optimality they have to intersect at the same value of $s$ since $s$ has to be minimum.  Using (\ref{rc2}),  we have 
$$
s \frac{\sinh(\sqrt{1-c_1^2} s)}{\sqrt{1-c_1^2} s}=
s \frac{\sinh(\sqrt{1-c_2^2} s)}{\sqrt{1-c_2^2} s}. 
$$ 
However the function $f(x)=\frac{\sinh(x)}{x}$ is strictly increasing with $x \in (0,\infty)$. Therefore, we must have $ {\sqrt{1-c_1^2} s}={\sqrt{1-c_2^2} s}$ which implies $c_2=\pm c_1$.

To find the point of crossing of the $x$ axis and study how it varies with $c$, let us fix $0< c <1  $, and consider the function for $y=y(s)$ given in (\ref{yt}) (\ref{cass1}), whose derivative is

$$
\frac{dy}{ds}=\frac{1}{\sqrt{1-c^2}}   \sinh (\sqrt{1-c^2} s)\sin(cs). 
$$

This derivative shows that $y$ grows until $s=\frac{\pi}{c}$ where it reaches its positive local maximum of $y_{max}=\frac{c}{\sqrt{1-c^2}}\sinh\left(\sqrt{1-c^2} \frac{\pi}{c} \right)$ and then decreases for $\frac{\pi}{c} < s< \frac{2\pi}{c}$ to reach at $s=\frac{2\pi}{c}$ its local minimum at $y_{min}=-\frac{c}{\sqrt{1-c^2}}\sinh \left( \sqrt{1-c^2} \frac{2\pi}{c} \right)$.
 The first intersection with the $x$ axis 
occurs for $\frac{\pi}{c} <  s  < \frac{3\pi}{2c}$ since $y(\frac{3\pi}{2c})=
-\cosh \left(\sqrt{1-c^2} \frac{3\pi}{2c}\right) <0$. The function $s_{int}:=s_{int}(c)$, giving the time where the intersection occurs,  is defined implicitly by $y(s_{int})=0$ 
for $\frac{\pi}{c} <  s=s_{int} < \frac{3\pi}{2c}$, that is, 
\be{sintc}
\tan(cs)=\frac{c}{\sqrt{1-c^2}}\tanh(\sqrt{1-c^2} s). 
\ee 
A detailed analysis based on the implicit function theorem given in Appendix B shows that $s_{int}=s_{int}(c)$ is a differentiable decreasing function of $c$, for $c \in (0,1)$ with $\lim_{c \rightarrow 0} s_{int}(c)=\infty$ and 
$\lim_{c \rightarrow 1}  s_{int}(c)=\bar s$, where $\bar s$ is the unique solution in 
$(\pi, \frac{3\pi}{2})$ of $\tan s=s$, which is,  $\bar s\approx  4.49341$. The coordinate  $x$ of the (first) intersection, $x_{int}$,   is negative and $|x_{int}|$ can be obtained from (\ref{rc2}), and it is 
$$
|x_{int}(c)|=\sqrt{1+\frac{\sinh^2(\sqrt{1-c^2} s_{int}{(c)})}{1-c^2}}. 
$$ 
This is, see Appendix B, a decreasing function of $c$, for $c \in (0,1)$ and, using the above limits for $s_{int}=s_{int}(c)$, we have 
$\lim_{c \rightarrow 0}|x_{int}(c)|=\infty$ and $\lim_{c \rightarrow 1} |x_{int}(c)|=\sqrt{1+\bar s^2}$, where $\bar s\approx  4.49341$ was described above. This  gives $\lim_{c \rightarrow 1} |x_{int}(c)|\approx 4.60333$. 

Figure \ref{Figura1} shows in green two optimal geodesics corresponding  to $0<c<1$ and their corresponding mirror image corresponding to $-c$. The geodesics on the left  correspond to $c=\pm 0.9$. The ones on the right correspond to $c=\pm 0.95$. The geodesics for $|c| \rightarrow 1$ converge to the ones for $c=\pm 1$ which are given in (\ref{xt}), (\ref{yt}), (\ref{traj2}) and plotted in black in the figure \ref{Figura1}.

\subsection{Geodesics for $|c|> 1$}

For $|c|>1$ , formula (\ref{rc2gen}) gives 

\be{rc2f}
r_c^2(s)=1+ \frac{1}{c^2-1} \sin^2(\sqrt{c^2-1} s). 
\ee
For a fixed $c$,  $r_c^2$ is equal to $1$ for the first time (after $s=0$) at $s=\frac{\pi}{\sqrt{c^2-1}}$. After this point the geodesic cannot be optimal anymore according to Proposition \ref{touch}, because it touches the singular part of the quotient space.  

As it was done for the case $|c|<1$,  we study the intersection (if any) of the geodesics (that are still symmetric with respect to the $x$ axis) with the $x$ axis. Because of symmetry, we consider only the case $c>0$. Similarly to what was done in the previous subsection, calculation of the derivative $\frac{dy}{ds}$ gives 
$$
\frac{dy}{ds}= \frac{1}{\sqrt{c^2 -1}}\sin(\sqrt{c^2-1} s) \sin(cs), 
$$ 
which shows that,  as long as $s< \frac{\pi}{c^2-1}$, $y$ is an increasing function until $s=\frac{\pi}{c}$ and   then decreasing until $s=\frac{2\pi}{c}$. The time $s=\frac{\pi}{c}$ always comes before the time $s=\frac{\pi}{\sqrt{c^2-1}}$ of intersection with 
the unit circle. However the time  $\frac{2\pi}{c}$ comes before $s=\frac{\pi}{\sqrt{c^2-1}}$ if and only if $c> \frac{2}{\sqrt{3}}$. The geodesic corresponding to $c=\frac{2}{\sqrt{3}}$ is such that at $s=\frac{\pi}{\sqrt{c^2-1}}=\frac{2\pi}{c}$, $r_c=1$ and $x=-1$. This gives a  pair of geodesics  in blue in figure \ref{Figura1} reaching the point $(-1,0)$.  For $1< c< \frac{2}{\sqrt{3}}$ the time  $s$ of intersection with the $x$ axis, $s_{int}=s_{int}(c)$ is at a value of $s$  which is  implicitly defined as a function of $c$ by setting $y=0$ in (\ref{yt}),  (\ref{cass2}). This function (which gives the sub-Riemannian distance of the points on the $x$ axis\footnote{This is identified with the minimum time in a time optimal control with bounded control.} ) is decreasing with $c$, for $ 1 < c< \frac{3}{2 \sqrt{2}}$ and then increasing with $c$, for $\frac{3}{2\sqrt{2}} < c \leq  \frac{2}{\sqrt{3}}$. 
We have continuity of $s_{int}(c)$ and $\lim_{c \rightarrow 1^+} s_{int}(c)=\bar s$, and 
$\lim_{c \rightarrow {\frac{2}{\sqrt{3}}}^-} s_{int}(c)=\sqrt{3} \pi$. 
The point of intersection $|x_{int}|=|x_{int}(c)|$ is decreasing as a function of $c$, with $\lim_{c \rightarrow 1^+}|x_{int}(c)|= \sqrt{1+\bar s^2}$ and $\lim_{c \rightarrow {\frac{2}{\sqrt{3}}^-}} |x_{int}(c)|=1$.  Detailed calculations justifying such statements are given in Appendix B. The geodesic corresponding to $c=\pm \frac{3}{2\sqrt{2}}$ is particularly interesting. It is plotted in red in figure \ref{Figura1}. Such geodesic  intersect the $x$-axis orthogonally at $s=s_{int}=\pi\sqrt{2}$, and a direct verification shows that  for $c=\pm \frac{3}{2\sqrt{2}}$ and 
for every $T$, $0 \leq T \leq \sqrt{2} \pi$,   
%
$$
(x_c,y_c)\left( \sqrt{2} \pi -T \right)=(x_{-c},y_{-c} ) \left( \sqrt{2} \pi + T \right). 
$$ 
Geodesics corresponding to $c> \frac{3}{2\sqrt{2}}$ curve inward and have an acute intersection angle. The figure reports in blue the geodesic corresponding to $c=\pm \frac{2}{\sqrt{3}}$, which includes the point $(-1,0)$. Also in blue are reported (on the left of the geodesic in red) the geodesic corresponding to $c=\pm 1.03$ which is such that $1 < |c| < \frac{3}{2\sqrt{2}}$  and  (on the right  of the geodesic in red) the geodesic corresponding to $c=\pm 1.12$ which is such that $\frac{3}{2\sqrt{2}}< |c| < \frac{2}{\sqrt{3}}$.

The geodesics corresponding to $|c| > \frac{2}{\sqrt{3}}$ are the ones that lose optimality by `landing' on the unit circle according to Proposition \ref{touch}.  The time of landing, which is the length of the sub-Riemannian geodesic $s_{int}$ , is given by the condition, obtained from (\ref{rc2f}),  $\sqrt{c^2-1} s_{int}=\pi$. Two such geodesics corresponding to $c=\pm  1.2$ and $c=\pm 1.5$ are depicted in purple  in the figure. The more the geodesics correspond to higher values of $|c|$ the shorter they are and the more their terminal point corresponds to a point on the circle with polar angle closer to zero. To see this, notice that for $t=\frac{\pi}{\sqrt{c^2-1}}$ the point (on the unit circle) is $(x_c,y_c)$, with $x_c=-\cos\left( \frac{c}{\sqrt{c^2-1}}\pi \right)$, $y_c=-\sin\left( \frac{c}{\sqrt{c^2-1}}\pi\right)$. As $c$ varies between $\frac{2}{\sqrt{3}}$ and $\infty$,  $\frac{c}{\sqrt{c^2-1}}$ decreases  monotonically from  $2$ to  $1$, and $x_c$ increases monotonically from $-1$ to $1$.

\subsection{Cut locus and Critical Locus}

The geodesics  do not intersect before reaching the $x$-axis or the unit circle.  More specifically, we have the 
following fact (cf. Appendix C).   

\bp{intersect} Denote by $s_{int}(c)$ the points where the geodesic corresponding to the value $c$ intersect the $x$-axis or the unit circle.   
Consider two optimal  geodesics $\gamma_{c_1}=\gamma_{c_1}(s)$ defined for $0 \leq s \leq s_{int}(c_1)$ and 
 $\gamma_{c_2}=\gamma_{c_2}(s)$ defined for $0 \leq s \leq s_{int}(c_2)$. They can only intersect when 
 $c_2=-c_1$ and at the $x$ axis with  $s_{int}(c_1)=s_{int}(c_2)=s_{int}(-c_1)$. 
\ep

Geodesics that reach the unit circle lose optimality on the unit circle according to Proposition \ref{touch}. Geodesics also lose optimality when intersecting the $x$-axis in points $(x,0)$ with $- \infty < x \leq -1$. In fact if the 
geodesic $\gamma_c$ did not lose optimality at the $x$-axis, the curve given by $\gamma_{-c}$ until 
the $x$-axis and $\gamma_{c}$ for an $\epsilon$-interval would still be optimal, contradicting smoothness. Lifting back to sub-Riemannian geodesics according to (\ref{critpoi}), we have the following proposition. 

\bp{critical}
The points where sub-Riemannian geodesics lose optimality in $SL(2)$ are in the inverse image under the natural projection of the union of the unit circle and the infinite interval $(-\infty, -1]$ on the $x$-axis.  
\ep

This is the set of the matrices in $SL(2)$ of the form (\ref{singpart7}) together with the matrices of the form (cf. the last paragraph of section \ref{Desk}) 
\be{newforma3}
X:\begin{pmatrix} a & b \cr b & d \end{pmatrix}, \qquad \texttt{with} \quad x=\frac{a+d}{2} \leq -1 .
\ee

For the sub-Riemannian problem a smoothness argument (see, e.g., \cite{AD}) shows that the 
sub-Riemannian cut locus ($SR-{\cal C L}$) is included in the sub-Riemannian critical locus. Therefore 
from (\ref{critpoi}) we have
$$
SR-{\cal CL} \subseteq \pi^{-1} \left(  R-{\cal CR}\right).   
$$

Matrices of the form (\ref{newforma3}) definitely belong to the cut locus since they are such that for the minimum $t_{min}$ they satisfy,  using (\ref{SRgeod}), 
$$
\pi \left( e^{(cA_0+P_1)t_{min}} e^{-cA_0 t_{min}} \right)= 
\pi \left( e^{(-cA_0+P_2)t_{min}} e^{cA_0 t_{min}} \right), 
$$ 
for matrices $P_1$ and $P_2$ in $\texttt{span} \{ A_1, A_2\}$. Therefore, there exists $K \in SO(2)$ such that
$$
e^{(-cA_0+P_2)t_{min}} e^{cA_0 t_{min}}= K e^{(cA_0+P_1)t_{min}} e^{-cA_0 t_{min}} K^T= e^{(cA_0+K P_1K^T )t_{min}} e^{-cA_0 t_{min}}. 
$$  
Therefore the two different  sub-Riemannian geodesics 
$\{ e^{(cA_0+K P_1K^T )t} e^{-cA_0 t} \, | \, t \in [0, t_{min}] \}$ and 
$\{ e^{(-cA_0+P_2)t} e^{cA_0 t} \, | \,  t \in [0, t_{min}] \}$ are both optimal and lead to the same point which therefore belongs to the sub-Riemannian cut locus.\footnote{Notice that no reference has been done to the specific form (\ref{newforma3}) and the same argument would show more in general that the inverse image of the Riemannian cut locus belongs to the sub Riemannian cut locus.}

Matrices of the form (\ref{singpart7}) also belong to the sub-Riemannian  cut locus because these 
matrices are obtained at the terminal time $t_{int}:=2s_{int}=\frac{2\pi}{\sqrt{c^2-1}}$ for the trajectories (\ref{SRgeod}). In particular,  a direct verification shows that the first factor in (\ref{SRgeod}) at $t_{int}$ is equal to $-{\bf 1} $,  independently of the matrix $P$   as long as $\|P\|=1$ (recall that $P$ has unit norm in (\ref{SRgeod})). Therefore for 
$\frac{2}{\sqrt{3}} < |c| < \infty$ the matrices 
$$
-e^{-cA_0 \frac{2\pi}{\sqrt{c^2-1}}}, 
$$
can all be obtained for an {\it infinite number} of geodesics each corresponding to  a matrix $P$ in 
the span of  
$A_1$ and $A_2$ and with $\|P\|=1$. 

In conclusion, we have 
\bp{CrCL}
For the sub-Riemannian structure of $SL(2)$ the critical locus and the cut locus 
coincide and are given by the inverse image under natural projection $\pi$ of the union of the  unit disc and the  semi-infinite segment $y\equiv 0$, $-\infty < x <-1$. These are the matrices described in (\ref{singpart7}), (\ref{newforma3}). 
\ep  

This result gives an alternative  description of the cut locus which is equivalent to the one given in \cite{BoscaRossi} (Theorem 5), with a different proof.\footnote{The matrices $K_{Id}^{loc}$ in \cite{BoscaRossi} correspond to the matrices (\ref{singpart7}) while the matrices $K_{Id}^{sym}$ correspond to (\ref{newforma3}). The proof that symmetric matrices in $SL(2)$ with trace $<0$ (as in \cite{BoscaRossi}) correspond to symmetric matrices with trace $<-2$ (as in (\ref{newforma3}))  boils down to showing that,  for positive $a$ and $d$, $ad>1$ implies $a+d>2$. In fact since $(a-1)^2 \geq 0$, $a +\frac{1}{a} \geq 2$ and since from $ad>1$, $d> \frac{1}{a}$ we have $a+d>a+\frac{1}{a}\geq 2$.}


\section{Finding the sub-Riemannian geodesics; An example}\label{Exe}

From a practical perspective, the main advantage of the characterization of geodesics described in the previous section is that it makes it easy to explicitly determine the sub-Riemannian geodesic between any  two points in $SL(2)$, that is, the parameter $c$ the matrix $P$ and the optimal final time $t$ in (\ref{SRgeod}). The fact that the problem was reduced to a two dimensional problem via symmetry reduction allows for a direct graphical  solution. 
The procedure, which is adapted from what described in \cite{AD},  can be summarized as follow: Given 
the desired initial and final conditions $X_i$ and $X_f$ in $SL(2)$, because of right invariance, the sub-Riemannian geodesic from  $X_{i}$  to  $X_f$ is $\gamma X_i$ where $\gamma$ is the sub-Riemannian geodesic from the identity to $\hat X_f:=X_f X_{i}^{-1}$. Therefore the problem is to determine $\gamma$. The algorithm is as follows: 

\begin{enumerate}
\item Determine the equivalence class of $\hat X_{f}$ ad therefore the final point $P_f$  in the plane of figure \ref{Figura1}. 

\item Determine the rough range of the parameter $|c|$ according to the desired final point, i.e., whether $|c|=1,0$  or $0 < |c| < 1$, or $1< |c| \leq \frac{2}{\sqrt{3}}$, or $|c|> \frac{2}{\sqrt{3}}$. In the various cases, the geodesic will be, respectively, from figure \ref{Figura1}, plotted in black, following the $x$ axis outwards starting from the point $(1,0)$, plotted in green, plotted in blue except for the geodesics corresponding to $|c|=\frac{3}{2 \sqrt{2}}$ which is plotted in red, 
plotted in purple.  

\item With a procedure of trial and error by considering geodesics corresponding to various values of $c$ in the given range (possibly with a bisection algorithm, see example below) find the value of $c$ which identifies the Riemannian geodesic crossing the point $P_f$. 

\item From the Riemannian geodesic,  determine the final optimal value $t_f$ where the point $P_f$ is reached. This is the sub-Riemannian distance. 

\item With formula (\ref{SRgeod}) and the found values of $c$ and $t_f$ using (for instance) $P=A_2$ find a given final condition, $Y_f$. This is such that $\pi(Y_f)=\pi(\hat X_f)$. 

\item Choose $K \in SO(2)$ so that $KY_f K^T=\hat X_f$. Then the correct $P$ is $P=KA_2K^T$.

\end{enumerate}

\subsection*{Example} 
\begin{figure}[h]
\centerline{\includegraphics[width=4.8in,height=3.5in]{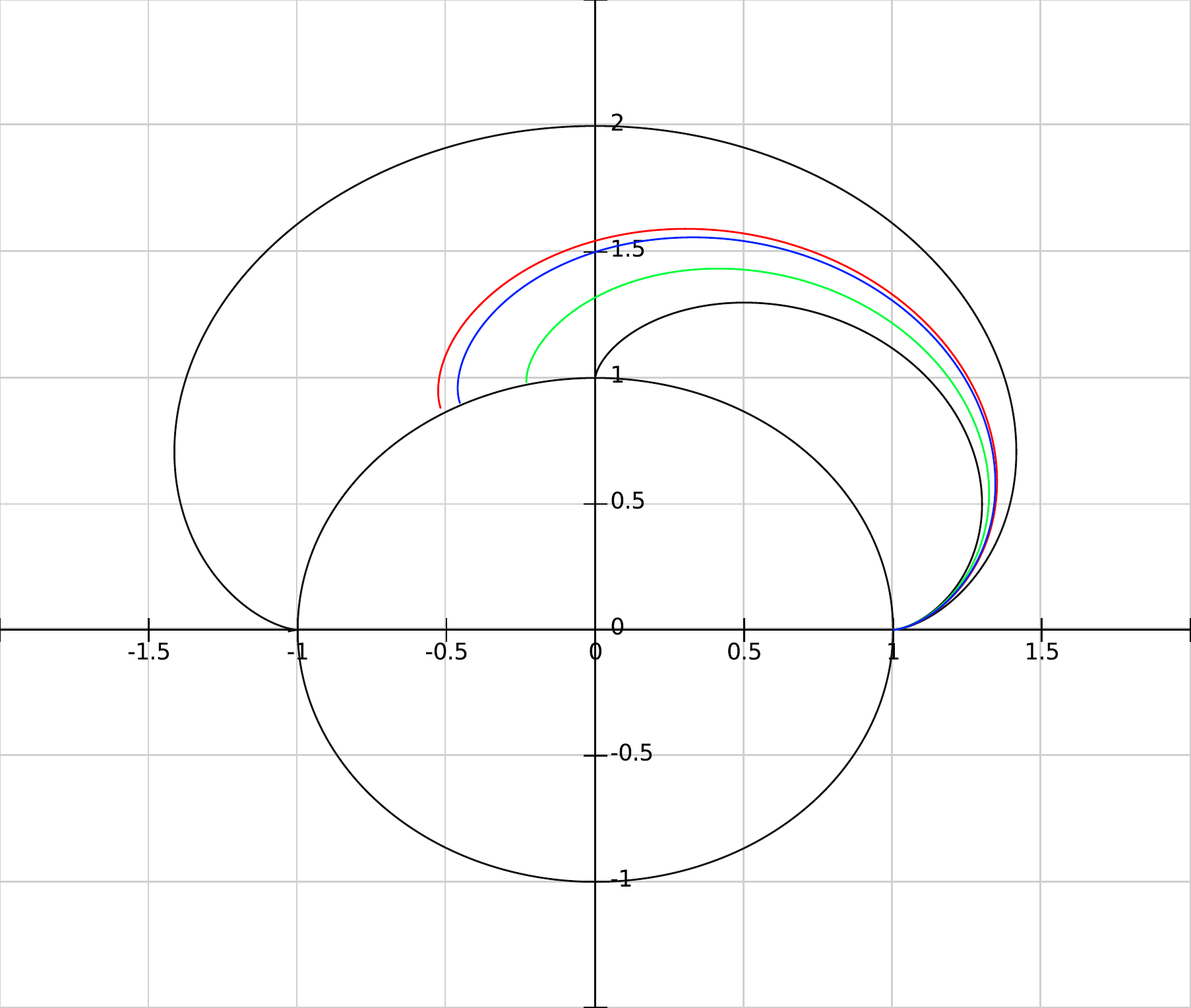}}
\caption{Riemannian geodesics for the example }
\label{Figura2}
\end{figure}

Assume we want to find the optimal sub-Riemannian geodesic from the point 
 $X_{i}=\begin{pmatrix} 0 & -1 \cr 1 & 0 \end{pmatrix} \in SL(2)$ to the point 
 $X_f=\begin{pmatrix} 2 & 1 \cr 1 & 1  \end{pmatrix} \in SL(2)$. Using (\ref{SRgeod}),  the geodesic has the form 
 $$
 X(t)=e^{(cA_0 +P)t}e^{-cA_0 t}X_i, \qquad t \in [0,t_f], 
 $$
and the problem is to find the minimum final time $t_f$, the optimal parameters $c$ and $P \in \texttt{span}\{A_1, A_2\}$, with $\|P\|=1$ such that $e^{(cA_0 +P)t_f}e^{-cA_0 t_f}=X_f X_i^{-1}$. In this case, we have 
\be{classe9}
X_f X_i^{-1}=\begin{pmatrix}- 1 & 2 \cr -1 & 1\end{pmatrix},  
\ee
which belongs to the class in $SL(2)/SO(2)$ corresponding to the point $(x,y)=(0,\frac{3}{2})$. By 
locating  this point on the diagram of Figure \ref{Figura1}, we find that the corresponding value 
of $c$ has to be positive and such that $\frac{2}{\sqrt{3}}< c < \infty$. To find the correct value of $c$ we plot different geodesics in the interval 
$\left( \frac{2}{\sqrt{3}}, \infty \right)$ looking for the one crossing the $y$ axis at the point with  $y=\frac{3}{2}$. In this case, we can in fact restrict the range of values of $c$ 
by excluding the geodesics which reach the unit circle (where they lose optimality)  before crossing the $y$ axis. In particular using (\ref{rc2f}) we know that a geodesic (\ref{xt}), (\ref{yt}), (\ref{cass2}),  reaches the unit circle for the first time at 
$s=\frac{\pi}{\sqrt{c^2-1}}$. Imposing that this happens when $x=0$ (using (\ref{xt}) (\ref{cass2})) gives that $c$ must satisfy $c=\frac{m}{\sqrt{r^2-4}}$ for an odd integer $m$. The only suitable odd integer which gives a value of $c$  
bigger than $\frac{2}{\sqrt{3}}$ is $m=3$ which gives $c=\frac{3}{\sqrt{5}}$. Therefore we consider the (Riemannian) geodesic from the point $(1,0)$ to the point $(0,1)$, and the sought after geodesic must be between this one and the one corresponding to $c=\frac{2}{\sqrt{3}}$ which is the one reaching the point $(-1,0)$. These two geodesics are in black in figure \ref{Figura2} and the value of $c$ must be such that $\frac{2}{\sqrt{3}}< c < \frac{3}{\sqrt{5}}$. 

Following step 3 of the above algorithm we first plot the Riemannian geodesic for $c=\frac{1}{2} \left( \frac{2}{\sqrt{3}} + \frac{3}{\sqrt{5}} \right)\approx 1.248171$, which gives the geodesic in red in Figure 2. Since this curve crosses the 
$y$ axis {\it above} the point $\frac{3}{2}$, the next value of $c$ is chosen (bisection) as $\frac{1}{2} \left( 1.248171 + \frac{3}{\sqrt{5}}\right)\approx 1.294906$ which gives the geodesic in green in the figure. The next approximation will be a curve between the red curve and green curve in the figure corresponding to a value of $c$, $1.248171 < c < 1.294906$. Proceeding this way we converge to the desired value of $c$ which is $c \approx 1.257558:=c_o$ which is the curve in blue in Figure $\ref{Figura2}$. 
 
The optimal value of $t$ (or $t_f$) can be found by setting $s=\frac{t}{2}$ and using (\ref{rc2f}) 
with $c=c_0$ and setting $r_c^2=\left( \frac{3}{2} \right)^2$. This gives for $s$ the value $s_0\approx 2.78115$. With these values we calculate (\ref{SRgeod}) with $t=2 s_0$ and $c=c_0$, and $P=A_2$ (point 5 of the algorithm) which gives 
$$
e^{(c_0 A_0+P)2 s_0} e^{-2cA_0s_0}=
$$
$$
\begin{pmatrix} \frac{1}{\sqrt{c_0^2-1}} \sin(\sqrt{c_0^2-1} s_0)+ \cos(\sqrt{c_0^2-1}s_0) &  -\frac{c_0}{\sqrt{c_0^2-1}} \sin(\sqrt{c_0^2-1} s_0)  \cr 
\frac{c_0}{\sqrt{c_0^2-1}} \sin(\sqrt{c_0^2-1} s_0) & - \frac{1}{\sqrt{c_0^2-1}} \sin(\sqrt{c_0^2-1} s_0)+ \cos(\sqrt{c_0^2-1}s_0) \end{pmatrix} \begin{pmatrix} \cos(c_0 s_0) & \sin(c_0 s_0) \cr - \sin(c_0 s_0) & \cos(c_0 s_0) \end{pmatrix} \approx  
$$
$$ 
\begin{pmatrix} 0.59540692401 & - 1.40601714179 \cr  1.40601714179  & 
-1.64070010186  \end{pmatrix} \begin{pmatrix}   -0.93732364113 & - 0.34846002895 \cr  0.34846002895  & -0.93732364113  \end{pmatrix} \approx 
$$
$$
\begin{pmatrix}-  1. 04 802976 & 1.11041756 \cr  
- 1.8896115 & 1.04792621 \end{pmatrix}:=Y_f
$$
This matrix is in the same equivalence class as (\ref{classe9}) but not exactly equal to it (although very close). According to step 6 of the algorithm, we need to find a matrix $K \in SO(2)$ such that 
$$
KY_fK^T=\begin{pmatrix} -1 & 2 \cr -1 & 1 \end{pmatrix}.
$$
This gives 
$$
K \approx \begin{pmatrix} 0.9170700563 & 0.39872611144 \cr 
-0.39872611144 & 0.9170700563
\end{pmatrix}.
$$
Therefore the correct matrix $P$ is 
$$
P=KA_2K^T \approx \frac{1}{2} \begin{pmatrix}  0.682034976 &  -0.73131955488  \cr -0.73131955488  & - 0.682034976  \end{pmatrix} 
$$


\section{Relation with the sub-Riemannian problem on $SU(2)$}\label{Relat5}
\begin{figure}[h]
\centerline{\includegraphics[width=4.8in,height=3.5in]{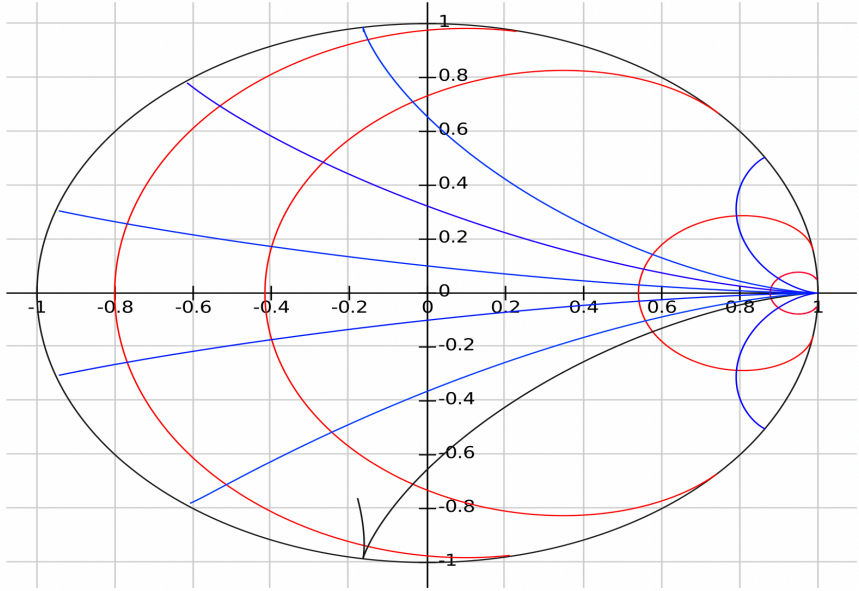}}
\caption{Riemannian geodesics in the $SU(2)$ case which are in blue. In red are depicted the boundaries of the reachable sets for various values of the time $s$. These are obtained by keeping $s$ fixed and varying $\omega$ in an appropriate range. We refer to \cite{AD}, \cite{S} for a discussion of the role of reachable set in optimal control problems in the context of symmetry reduction.}
\label{Figura3}
\end{figure}

In \cite{Adpre} the sub-Riemannian $K-P$ problem for the Lie group $SU(2)$ was solved using the technique of symmetry reduction advocated here. It is interesting to compare the results. The Lie algebras $su(2)$ and $sl(2)$ are related in that they are both real forms of the complex Lie algebra $sl(2,C)$.\footnote{Recall that for a complex Lie algebra ${\cal L}$,  we can consider ${\cal L}$ as a Lie algebra over the reals which has twice the (complex) dimension as ${\cal L}$ and we denote it by ${\cal L}^R$. The real Lie algebra ${\cal L}^R$ has a {\it complex structure} $J$, that is, an endomorphism ${\cal L}^R \rightarrow {\cal L}^R$ such that $J^2=-Id$, with $Id$ denoting the identity, given by multiplication  by the imaginary unit $i$. A real form of ${\cal L}$ is a real Lie algebra ${\cal L}_0$ such that ${\cal L}^R={\cal L}_0 \oplus J {\cal L}_0$. Both $sl(2)$ and $su(2)$ play this role for $sl(2,C)$. } The Lie algebra $su(2)$ also has a $K-P$ Cartan decomposition, $su(2)={\cal K} \oplus {\cal P}$  with 
${\cal K}=\texttt{span}\{A_0\}$ and ${\cal P}:=\texttt{span}\{ i A_1, iA_2 \}$ (cf. (\ref{f0f1f2})). The metric for the sub-Riemannian problem is defined as in (\ref{metrica}) with transposition $^T$ replaced by transpose conjugate $^\dagger$. In both cases the action of the Lie group $SO(2)$ can be taken as the symmetry.   The invariants  (parametrizing the quotient space) under the (conjugation) action $X \rightarrow KXK^T$ ($K \in SO(2)$) are defined in the same way as in (\ref{genmatX}) (\ref{xey}),  $x:=\frac{Tr(X)}{2}$, $y=\frac{X_{12}-X_{21}}{2}$,\footnote{Notice that every matrix in $SU(2)$ can be written as $X=\begin{pmatrix} w & z \cr - z^* & w^* \end{pmatrix}$, therefore both $x$ and $y$ are real.} which are in this case constrained {\it inside} the unit disc (unit circle included) in the $x-y$ plane \cite{Adpre}.  The quotient map $\pi$ has the same form as $\pi$ in (\ref{pi}) but now it acts on $SU(2)$. The sub-Riemannian geodesics have the form (\ref{SRgeod}) but with $P$ replaced by $-iP$. Therefore, the sub-Riemannian geodesics  in the $SU(2)$ case can be obtained from (\ref{SRgeod}) by the substitution $c \rightarrow i \omega$, $t \rightarrow -is$ and are therefore given by $X^U(s):=e^{(A_0 -iP)s}e^{-\omega A_0s}$, which are parametrized by a parameter $\omega$.  Since the projection $\pi$ takes the same form, Riemannian geodesics in the $SU(2)$ case\footnote{From  the analysis in \cite{DS}  for the $SU(2)$ case the regular part in $SU(2)/SO(2)$ corresponds to the {\it interior} of the unit disc in the plane.}  
can be obtained by doing the same substitution $c \rightarrow i \omega$, $ t \rightarrow -is$ in the   Riemannian geodesic for the $SL(2)$ case. To this purpose, we notice that formulas (\ref{xt}) (\ref{yt}) with the choices (\ref{cass1}) (\ref{cass2})  are both defined when $c=i\omega$ independently of the absolute value of $|c|=|\omega|$ and give the same expressions for the Riemannian geodesics in the $SU(2)$ case. So, 
for example, we can use the expression of the  Riemannian geodesic for $|c| > 1$ in the $SL(2)$ case (obtained from (\ref{xt}), (\ref{yt}), (\ref{cass2}) replacing the variable $s$ with $t$, with some abuse of notation)  
\be{geoL}
\gamma_c^L:= \left\{ \begin{matrix} x=\cos(\sqrt{c^2 -1} t) \cos(ct)+ \frac{c}{\sqrt{c^2-1}} \sin(\sqrt{c^2-1}t) 
\sin(ct)\\
y=\cos(\sqrt{c^2-1} t) \sin(ct) - \frac{c}{\sqrt{c^2 -1}} \sin(\sqrt{c^2-1}t) \cos(ct) \end{matrix} \right., 
\ee
and we obtain the Riemannian geodesic in the $SU(2)$ case, 
\be{geoU}
\gamma_\omega^U=\left\{ \begin{matrix} x=\cos(\sqrt{1+ \omega^2} s) \cos(\omega s) + \frac{\omega}{\sqrt{\omega^2+1}} \sin(\sqrt{\omega^2+ 1}s) \sin(\omega s)  \\ 
y= \cos(\sqrt{1+ \omega^2} s) \sin(\omega s) - \frac{\omega}{\sqrt{\omega^2+1}} \sin ( \sqrt{\omega^2+1} s) \cos(\omega s) \end{matrix} \right. 
\ee
for all values of $\omega$. For a fixed value of $\omega$ the Riemannian geodesic in the $SU(2)$ case 
starts at the $(1,0)$ point and ends, losing optimality,  on  the unit circle, which, similarly to the $SL(2)$ case, is the singular part of the quotient space.  A plot of the Riemannian geodesics in the $SU(2)$ case in figure \ref{Figura3}, reveals that there is a one to one correspondence  between geodesics and values of $\omega$ and all values of $\omega$ are possible with a symmetry $\omega \leftrightarrow -\omega$ about the $x$-axis. The geodesic  corresponding  to $\omega=0$ connects going leftwards the point $(1,0)$ with the point $(-1,0)$, while the others, as $|\omega|$ increases,  land on a point of the unit circle whose polar coordinate angle is closer to zero, losing optimality.

It is possible to obtain {\it all} the $SU(2)$ Riemannian geodesics (all possible values of $\omega$) as corresponding to $SL(2)$ Riemannian geodesics with  a range of values of $c$ different from $(-\infty, \infty)$. Furthermore, we can establish a correspondence between $SU(2)$ and $SL(2)$ Riemannian geodesics so that the corresponding geodesics  land on the same point on the unit circle.  In order to do  this,  assume first $\omega \geq 0$. Then the corresponding value of $c$ is chosen as 
\be{choicc}
c=-\sqrt{ \frac{5 \omega^2 +4 - 4 \omega \sqrt{\omega^2+1}}{4 \omega^2+3 -4 \omega \sqrt{\omega^2+1}}}.  
\ee 
This is a monotonically decreasing finction of $\omega$ from $c=-\frac{2}{\sqrt{3}}$ for $\omega=0$ to $-\infty$ for $\omega \rightarrow \infty$. Formula (\ref{choicc}) after squaring both terms and direct manipulations gives 
$$
\frac{c^2}{c^2-1}=\frac{(\omega-2 \sqrt{\omega^2+1})^2}{\omega^2+1}, 
$$ 
which gives 
\be{true2}
\frac{c\pi}{\sqrt{c^2-1}}=-2\pi + \frac{\omega \pi}{\sqrt{\omega^2 +1}}. 
\ee
The point where the $SU(2)$  geodesic (\ref{geoU}) reaches the unit circle is the point 
$$
x=-\cos\left( \frac{\omega \pi}{\sqrt{\omega^2+1}}\right), \qquad 
y=-\sin\left( \frac{\omega \pi}{\sqrt{\omega^2+1}}\right). 
$$
The point where the $SL(2)$  geodesic (\ref{geoL}) reaches the unit circle  is the point 
$$
x=-\cos\left( \frac{c \pi}{\sqrt{c^2-1}}\right), \qquad 
y=-\sin\left( \frac{c \pi}{\sqrt{c^2-1}}\right). 
$$
Because of (\ref{true2}) the correspondence (\ref{choicc}) guarantees that these two points coincide. An analogous calculation for $\omega \leq 0$ shows that the choice 
$$
c=\sqrt{ \frac{(\omega +2 \sqrt{\omega^2+1})^2}{4\omega^2+3+4 \omega \sqrt{\omega^2+1}}}, 
$$
which is a monotonically decreasing function (as $\omega$ goes from  $0$ to  $ -\infty$,  $c$ goes from $\frac{2}{\sqrt{3}}$ to $+\infty$) ensures that the two geodesics corresponding to $\omega$ and $c$ end (and lose optimality)  at the same point of the unit disc.

\section{Automorphic structures on $SL(2)$}\label{Ident5}

As we have mentioned, sub-Riemannian structures $(M, \Delta, g)$, $(M, \Delta^{'}, g^{'}))$, with $M$ a Lie group,  can be considered as equivalent if there exists an automorphism $\phi$ of $M$, with $\phi_* \Delta= \Delta^{'}$  and $\phi^{*} g^{'}=g$. In this case, the structures can be identified and geodesics of 
one structure can be mapped to geodesics in the other structure via the automorphism $\phi$. It is therefore of interest to identify the structures which are related by an automorphism to the one we have studied in the previous sections,  since for these structures as well we will have  found the complete optimal synthesis. 

The study of the automorphisms of  Lie groups  is an important chapter  of the theory of Lie groups (see, e.g., \cite{Helgason}, and \cite{Dani} for a recent expository paper). In the following, we always assume that $M$ is compact. The group of automorphisms on $M$, $Aut(M)$,  can itself be given the structure of a Lie group.  A subgroup of $Aut(M)$ is given by the group of {\it inner automorphisms}, $Inn(M)$. These are automorphisms $\phi_I$ such that there exists an $m \in M$ so that $\phi_I(x)=mxm^{-1}$ for every $x\in M$. The subgroup $Inn(M)$ is a normal subgroup of $Aut(M)$, and the quotient ${Aut(M)}/{Inn(M)}$ is called the {\it group of outer automorphisms}. For a connected and simple Lie group, such as $SL(2)$, the group $Inn(M)$ is the connected component of $Aut(M)$ containing the identity. The map $\phi \rightarrow \phi_*$, applied at the identity in $M$,   maps  an automorphism $\phi \in Aut(M)$ to a Lie algebra automorphism $\phi_*$. The group of Lie algebra automorphisms $aut\left({\cal M}\right)$, with ${\cal M}$ the Lie algebra of ${M}$,  is the group of all invertible linear maps ${\cal M} \rightarrow {\cal M}$ that preserve the Lie bracket, i.e.,  $\alpha \in aut\left( {\cal M} \right) \leftrightarrow \alpha \left( [A,B] \right)= [ \alpha(A), \alpha(B)]$. The map $\phi \rightarrow \phi_*$ from $Aut(M)$ to $aut\left( {\cal M} \right)$ is injective but not necessarily surjective. It is surjective if $M$ is simply connected and it is also surjective, we will see, in the $SL(2)$ case, although $SL(2)$ is not simply connected. The automorphisms of all classical Lie groups have been described in detail in standard  references \cite{Dieud}.   Here we present,  for the $SL(2)$ case, an elementary and explicit derivation based on the above map $\phi \rightarrow \phi_*$. In particular, we first describe the group $aut \left( sl(2) \right)$ using the argument in \cite{HB2} and then we find $Aut \left( SL(2) \right)$ by displaying for every $\alpha \in   aut \left( sl(2) \right)$ a 
$\phi \in Aut \left( SL(2) \right)$ such that $\phi_*=\alpha$.

In describing $aut \left( sl(2) \right)$, we consider the coordinates of $sl(2)$ in the basis $\{ A_0 ,A_1,A_2\}$. With this basis,  the sub-Riemannian structure we have treated above is identified by the pair of vectors $\{ \vec e_2, \vec e_3\}$ in $\RR^3$. The automorphism group $aut \left( sl(2) \right)$ is identified with a subgroup of   the linear group on $\RR^3$. We have the following. 
\bl{autom}
$$
aut\left( sl(2) \right) = SO_0(1,2)
$$
\el  
The Lie group $SO(1,2)$ is  the group of $3 \times 3$ matrices $X$ preserving the quadratic form identified by the matrix 
$I_{1,2}:=\begin{pmatrix} -1 & 0 & 0 \cr 0 & 1 & 0 \cr 0 & 0 & 1 \end{pmatrix} $, i.e., $X^T I_{1,2} X=I_{1,2}$. This Lie group (and its extension $SO(1,3)$) is also known as the Lorentz group and has applications in relativity. $SO(1,2)$ has two connected components (Lemma X.2.4 in \cite{Helgason}) given by matrices with determinant equal to $+1$ or $-1$, respectively. The component $SO_0(1,2)$ corresponds the matrices with determinant equal to $1$. This is the component which contains the identity. 

Lemma  \ref{autom} was proven  in \cite{HB2} and we will give the proof with some 
extra details in Appendix D.

We now give a concrete parametrization of matrices in $SO_0(1,2)$ (cf. also, section X.2.3 in \cite{Helgason}).

\bl{Paramet}
Every matrix $X$ in $SO_0(1,2)$ can be written as 
\be{factoriz}
A= O(\theta_1) I^{0,1,2} H(z) O(\theta_2), 
\ee 
for three real parameters $\theta_1, \theta_2, z$ and one of three  choices $0,1,2$, where,  
\be{Oteta}
O(\theta):=\begin{pmatrix}  1 & 0 & 0 \cr 0 & \cos(\theta) & \sin(\theta) \cr 0 & -\sin(\theta) & \cos(\theta) \end{pmatrix}, \quad H(z)=\begin{pmatrix} \cosh(z) & 0 & \sinh(z) \cr 0 & 1 & 0 \cr \sinh(z) & 0 & \cosh(z)\end{pmatrix}, 
\ee
and  $I^0$ is the $3 \times 3$ identity,  $I^1:=\texttt{diag} (-1, -1, 1)$, $I^2:=\texttt{diag} (-1, 1,-1 )$.  
\el 
Notice that the definition of $I^{0,1,2}$ is somehow redundant, since 
$I^2=O\left( \frac{\pi}{2} \right) I^1 O\left( -\frac{\pi}{2} \right)$, but it was set this way to simplify the exposition in the proof below. 
\bpr 
Given a matrix $A \in SO_0(1,2)$ choose $\hat \theta_1$ so that $O(\hat \theta_1) A$ has the $(2,1)$ entry equal to zero and then choose $\hat \theta_2$ so that $O(\hat \theta_1) A O(\hat \theta_2)$ has the $(1,2)$ entry equal to zero (notice this does not affect $(1,2)$ entry). Therefore we have 
\be{RX2}
O(\hat \theta_1) A O(\hat \theta_2)=\begin{pmatrix} h_{11} & 0 & h_{13} \cr 0 & h_{22} & h_{23} \cr h_{31} & h_{32} & h_{33}  \end{pmatrix}:=\hat H. 
\ee
Since all the matrices on the left hand side of (\ref{RX2}) are in $SO_0(1,2)$, so is the matrix $\hat H$, which therefore satisfies $\hat H^T I_{1,2} \hat H=I_{1,2}$. From this and (\ref{RX2}), we obtain the relations 
\begin{eqnarray}
-h_{11}^2 + h_{31}^2= -1 \label{Euno} \\
h_{31} h_{32}=0 \label{Edue}\\
-h_{11}h_{13}+ h_{31} h_{33}=0 \label{Etre} \\
h_{22}^2+h_{32}^2 =1 \label{Equattro}\\
h_{22} h_{23} + h_{32} h_{33} =0 \label{Ecinque}\\
- h_{13}^2 + h_{23}^2 + h_{33}^2=1. \nonumber
\end{eqnarray}
From (\ref{Edue}) we know that at least one between $h_{31}$ and $h_{32}$ is zero. If $h_{31}=0$, from (\ref{Etre}) and (\ref{Euno}),  we have that $h_{13}=0$. In that case,  $\hat H$ has the form $I^{0,1} O(\hat \theta_3)$ for some real parameter $\hat \theta_3$. If $h_{32}=0$, (\ref{Ecinque}) and (\ref{Equattro}) give $h_{23}=0$, and $h_{22}=\pm 1$. Thus, the matrix $\hat H$ reduces to 
$$
\hat H=\begin{pmatrix} h_{11} & 0 & h_{13} \cr 0 & \pm 1 & 0 \cr h_{31} & 0 & h_{33} \end{pmatrix}, 
$$
with 
\be{ridot}
-h_{11}^2 + h_{31}^2=-1, \qquad h_{31} h_{33}=h_{11} h_{13}, \quad h_{33}^2 -h_{13}^2=1. 
\ee
Squaring the second equation of (\ref{ridot}) and replacing the first and the third, we 
obtain $h_{31}^2=h_{13}^2$ which also gives $h_{11}^2=h_{33}^2$. Let us first consider the case $h_{13}=0$ and therefore $h_{31}=0$. We have $h_{11}=\pm 1$ and $h_{33}=\pm 1$. Taking into account that $\det (\hat H)=1$, there are only four  possibilities for $\hat H$: 1) $\hat H=\texttt{diag}(1,1,1)$, that is,  the identity; 2) $\hat H=\texttt{diag}(-1,-1,1)$, that is $I^1$; 3) $\hat H=\texttt{diag}(1,-1,-1)$, that is, the matrix $O(\pi)$, with the definition (\ref{Oteta}); 4) $\hat H=\texttt{diag}(-1,1,-1)$, that is, $\hat H=I^2$. 

If $h_{13}\not=0$, setting $h_{13}= \sinh( z)$, for some real $z$,  we consider two possibilities $h_{11}>0$ and $h_{11}<0$. If $h_{11}>0$, from the first one of (\ref{ridot}) and $h_{31}^2=h_{13}^2$ we can write $h_{11}=\cosh(z)$ and $h_{33}=\pm \cosh(z)$, $h_{31}=\pm \sinh(z)$. From the second one of (\ref{ridot}) it follows that the same sign has to be used in $h_{31}$ and $h_{33}$. The condition $\det(\hat H)=1$ imposes that if the chosen sign is $+$ then $h_{22}=1$ and $\hat H$ has exactly 
the form of $H$ in (\ref{Oteta}). If we choose the sign $-$, then, the condition   $\det(\hat H)=1$ imposes that $h_{22}=-1$ and again $\hat H=O(\pi) H$ where $O$ and $H$ were defined in (\ref{Oteta}).  In the case,  $h_{11}<0$,  we write $h_{11}=-\cosh(z)$, and again $h_{33}=\pm \cosh(z)$, 
$h_{31}=\pm \sinh(z)$. With these definitions, the second formula of (\ref{ridot}) gives that the opposite signs have to be chosen. By choosing $-$ for the $\cosh$ and $+$ for the $\sinh$, the condition on the determinant gives that $h_{22}=+1$ and therefore $\hat H=I^2 H(-z)$. By choosing $+$ for the $\cosh$ and $-$ for the $\sinh$ we have that $h_{22}=-1$ and  $\hat H=I^1 H(-z)$. In all cases considered above for $\hat H$, from formula (\ref{RX2}) we have that $A$ has the factorization (\ref{factoriz}) (with 
$\hat \theta_{1,2}=-\theta_{1,2})$. 
\epr 

The following proposition describes the group of automorphisms on $SL(2)$, that is, $Aut \left( SL(2) \right)$. 

\bp{Autom90}
The group of automorphisms on $SL(2)$, $Aut \left( SL(2) \right)$, is given by all the maps 
$\phi(X)=KXK^{-1}$ with $K$ belonging to $SL^{\pm}(2)$,  the Lie group of $2 \times 2$ real matrices with determinant equal to $\pm 1$. 
\ep  
Notice that $Aut \left( SL(2) \right)$ is different from $Inn \left( SL(2) \right)$ and the group of outer automorphisms, that is,  
\newline 
 ${Aut\left( SL(2) \right)}/{Inn\left( SL(2) \right)}$,  is a cyclic group of order two with two classes given by the class containing the identity (and corresponding to the inner automorphisms)  and the class containing the transformation $\phi_-(X)=K_- X K_-^{-1}$ where $K_-:=\begin{pmatrix} 0 & 1 \cr 1 & 0 \end{pmatrix} \in SL^{-}(2)$. 
\bpr 
The map $\phi \rightarrow \phi_*$ from $Aut(SL(2))$ to $aut(sl(2))$ is an injective homomorphism  and thus one can naturally consider $Aut(SL(2))$ as a subgroup of $aut(sl(2))$ (see, e.g., \cite{Dani}). The group $G$ defined in the statement of the theorem is a subgroup of $Aut(SL(2))$. Thus the proposition is proved if we prove that for every $\alpha \in aut(sl(2))$ there exists a $\phi \in G$ such that $\phi_*=\alpha$. We notice that the automorphisms $\phi$  in $G$ on $SL(2)$ have formally the same (conjugacy) expression if we consider the corresponding a $\phi_*$ on $sl(2)$, that is $A \rightarrow KAK^{-1}$. 
Therefore, the proposition is proved, using Lemma \ref{autom}, if we show that for every element  in $SO_0(1,2) \simeq aut\left( sl(2) \right)$, there exists a $K \in SL^{\pm}(2)$ which  realizes the same transformation  via $A \rightarrow KAK^{-1}$ $A \in sl(2)$. In fact,  it is enough to show this for the three types of factors   which appear in (\ref{factoriz}), $O(\theta)$, $H(z)$ and $I^{0,1,2}$, since,  according to Lemma \ref{Paramet},  they generate $SO_0(1,2)$. We do not need to check this for $I^0$ which is just the identity and for $I^1$ that can be obtained from $I^2$ and $O(\theta)$ as $O\left(\frac{\pi}{2}\right) I^2O\left(-\frac{\pi}{2}\right)$.

The automorphisms in $aut(sl(2))$ of the type $O(\theta)$ are obtained as $A \rightarrow KAK^{-1}$ with $K=\begin{pmatrix} \cos\left( \frac{\theta}{2} \right) & \sin\left( \frac{\theta}{2} \right) \cr 
- \sin\left( \frac{\theta}{2} \right)  & \cos\left( \frac{\theta}{2} \right) 
 \end{pmatrix}$, since, with this choice,  direct calculation shows $KA_0 K^{-1}=A_0$, $KA_1 K^{-1}= \cos(\theta) A_1+\sin(\theta) A_2$, $KA_2K^{-1}=-\sin(\theta) A_1 + \cos(\theta) A_2$ which, in coordinates,  is the transformation $O(\theta) \in SO_0(1,2)$. The automorphisms in $aut(sl(2))$ of the type $I^1$ are obtained with $K=\begin{pmatrix} 0 & 1 \cr 1 & 0 \end{pmatrix}$ since we have, with this choice,  $KA_0 K^{-1}=-A_0$, $KA_1 K^{-1}=A_1$, $KA_2 K^{-1}=-A_2$. Finally, the automorphism $H(z)$ in (\ref{Oteta}) is obtained using $K=\begin{pmatrix} \cosh\left( \frac{z}{2} \right) & \sinh\left( \frac{z}{2}\right)  \cr  \sinh\left( \frac{z}{2}\right) &   \cosh\left( \frac{z}{2} \right) \end{pmatrix}$, as a similar direct verification shows. 

\epr 
In the $SL(2)$ case, the group $Iso_{\bf 1} \left( SL(2), \Delta, g \right)$ is a subgroup of $Aut \left( SL(2)\right)$\cite{Biggs}. \footnote{Even if we restrict to  right invariant sub-Riemannian structures on {\it three dimensional}  Lie groups it is not always true that   $Iso_{\bf 1} \left( M, \Delta, g \right)$ is a subgroup of $Aut(M)$. See  the `exceptional case' in \cite{Biggs}.} The  group $Iso_{\bf 1} \left( SL(2), \Delta, g \right)$ can be obtained from the representation (\ref{factoriz}) and consists of all the matrices in $SO_0(1,2)$ of the form $I^{0,1} O(\theta)$ for free real $\theta$. 

The structures on which our results in the previous sections apply can be identified  by two elements in $sl(2)$ $\{ B_1, B_2\}$ and their associated right invariant vector fields. The metric is implicitly defined by imposing that these two matrices are orthonormal. Our results apply to all structures  $\{ B_1, B_2\}$ automorphic to the structure $\{A_1, A_2\}$, that is, such that there exists a $K \in SL^{\pm} (2)$ with $KA_1K^{-1}=B_1$ and $KA_2 K^{-1}=B_2$.  
Equivalently,  these are structures whose coordinates in the $\{A_0, A_1,A_2\}$ basis are given by the last two columns of a matrix in $SO_0(1,2)$. For example, choosing a matrix of the form $H(z)$ in (\ref{Oteta}), we might have $B_1=A_1$, $B_2=\sinh(z) A_0 + \cosh(z) A_2$, for some real $z$. Let $\{ \vec a_1, \vec a_2 \}$ ($\{ \vec b_1, \vec b_2 \}$)  be the coordinate vectors corresponding to $\{A_1, A_2\}$ ($\{B_1, B_2\}$). We have 
$\vec a_j^T I_{1,2}  \vec a_k= \delta_{j,k}$ and since $\vec b_{j,k}=X \vec a_{j,k}$ with $X \in SO_0(1,2)$      
$\vec b_j^T I_{1,2}  \vec b_k= \delta_{j,k}$. In particular the inner product 
$\langle \vec v, \vec w \rangle= \vec v^T I_{1,2} \vec w$ is positive definite on any pair of vectors representing a sub-Riemannian structure automorphic to the one we have treated. Structures where such inner product is positive definite are called {\it elliptic}, and therefore the structures considered in this paper are a special case of elliptic sub-Riemannian structures.  In fact,  the inner product  $\langle \vec v, \vec w \rangle= \vec v^T I_{1,2} \vec w$ is the same as the {\it Killing} inner product defined by $\langle  A, B \rangle=Tr(ad_A ad_B)$  
(cf. (\ref{addav})) which is equal, up to a scaling factor,  to the inner product defined in (\ref{metrica}). Structures on $SL(2)$ such that the Killing inner product is indefinite are called {\it hyperbolic} \cite{AGBD}. Sub-Riemannian structures that have the same geodesics as described here (or mapped to them by an automorphism) are also the ones corresponding to pairs of coordinate vectors 
$\{ \lambda  X\vec a_1, \lambda X \vec a_2 \}$, for some scaling factor $\lambda > 0$, since the scaling factor can be taken into account be rescaling the metric which corresponds  to a uniform scaling of the sub-Riemannian length on the whole manifold. The explicit description of geodesics for sub-Riemannian structures (elliptic or hyperbolic) non automorphic to the one considered here is an open problem,  although important general properties have been proved \cite{Mashtakov}.

 \section*{Acknowledgement} D. D'Alessandro's  research was supported by NSF under Grant EECS-1710558

\section*{Appendix A: Solutions of the geodesic equations} 

We verify here that the trajectories (\ref{xt}) (\ref{yt}) solve the differential equations (\ref{geodeq1}), i.e., the geodesic equations in the regular part of the quotient space. Equivalently, we verify equations (\ref{geodeq1}) multiplied by $x^2+y^2-1>0$. Using (\ref{xt}) and (\ref{yt}) we obtain 
\be{x2y2}
x^2+y^2-1=k_1^2+k_2^2-1. 
\ee 
We shall also need the derivatives of $x=x(t)$ and $y=y(t)$. Defining  $s:=\frac{t}{2}$, we have 
$
\dot x=\left(\dot k_1 +\frac{k_2c}{2} \right)\cos(cs)-\left(\dot k_2-\frac{k_1c}{2} \right) \sin(cs). 
$
Using (\ref{traj1}) and (\ref{traj2}), it is straightforward to verify that in all cases $\dot k_2-\frac{k_1c}{2} =0$. Therefore we have 
\be{dotx}
\dot x=\left(\dot k_1 +\frac{k_2c}{2} \right)\cos(cs). 
\ee
Analogously, we find, 
\be{doty}
\dot y=\left(\dot k_1 +\frac{k_2c}{2} \right)\sin(cs). 
\ee
As for the second derivatives,  we get 
$
\ddot x=\left( \ddot k_1+ \frac{\dot k_2 c}{2} \right) \cos(cs) - \left( \dot k_1 +\frac{k_2 c}{2} \right) \frac{c}{2} \sin (cs), 
$ and 
$
\ddot y= \left( \ddot k_1+ \frac{\dot k_2 c}{2} \right) \sin(cs) + \left( \dot k_1 +\frac{k_2 c}{2} \right) \frac{c}{2} \cos (cs). 
$
We can verify with the explicit expressions of $k_1$ that $ \ddot k_1+ \frac{\dot k_2 c}{2} =\frac{k_1}{4}$ and 
$\left( \dot k_1 +\frac{k_2 c}{2} \right) \frac{c}{2}=\frac{k_2}{4}$. Therefore we have 
$$
\ddot x=\frac{k_1}{4} \cos(cs)-\frac{k_2}{4} \sin(cs), 
$$
$$
\ddot y=\frac{k_1}{4} \cos(cs)+\frac{k_2}{4} \sin(cs).
$$
Let us now verify the first one of (\ref{geodeq1}). The verification of the second one is analogous. Using (\ref{x2y2}) the left hand side of (\ref{geodeq1}) is 
\be{LHS}
(LHS)_x=\frac{k_1^2+k_2^2-1}{4}\left(k_1 \cos(cs)-k_2 \sin(cs) \right). 
\ee  
Using (\ref{xt}), (\ref{yt}), (\ref{dotx}), (\ref{doty}), the right hand side is 
\be{RHS}
(RHS)_x= \left[\dot k_1 +\frac{k_2c}{2} \right]^2 \left[ \left( k_1 \cos(cs)+k_2 \sin(cs)\right)\left( \cos^2(cs) -\sin^2(cs) \right)+2 \left(k_1 \sin(cs)-k_2 \cos(cs) \right) \cos(cs) \sin(cs)  \right]
\ee
$$
 =\left[\dot k_1 +\frac{k_2c}{2} \right]^2 \left[ k_1 \cos(cs) -k_2 \sin(cs) \right].
$$
Therefore, comparing (\ref{RHS}) with (\ref{LHS}), equality $(LHS)_x=(RHS)_x$ is proven if we prove  $\left[\dot k_1 +\frac{k_2c}{2} \right]^2=\frac{k_1^2+k_2^2-1}{4}$. However, this is verified directly using the expressions of $k_1$ and $k_2$ in (\ref{traj1}) and (\ref{traj2}). 

\section*{Appendix B: Analysis of the functions $s_{int}$ and $|x_{int}|$ }

\subsection*{Case $0<c<1$}

We shall use the following fact: 

\vs

{FACT:} {\it $\tanh(w)-w<0 $  for every $w>0$}. 

\vs

Define $k:=\frac{\sqrt{1-c^2}}{c}$, which is equivalent to $c=\frac{1}{\sqrt{1+k^2}}$,  and $x:=cs$ (with $\pi < x < \frac{3}{2} \pi$), so that 
\be{dkdc}
\frac{dk}{dc}=-\frac{1}{c^2 \sqrt{1-c^2}}=-\frac{(1+k^2)^{\frac{3}{2}}}{k} < 0. 
\ee

\subsubsection*{Function $s_{int}=s_{int}(c)$}

For simplicity we write $s(c)$ instead of $s_{int}(c)$. We have 
\be{derivata}
\frac{dx}{dc}=s+c\frac{ds}{dc}=\frac{dx}{dk}\frac{dk}{dc}=
-\frac{dx}{dk}\frac{(1+k^2)^{\frac{3}{2}}}{k}. 
\ee
The function $x=x(k)$ is defined implicitly by 
\be{Ffan}
F(k,x)=k\tan(x)-\tanh(kx)=0, 
\ee
for $\pi < x< \frac{3}{2} \pi$, which is obtained from  (\ref{sintc}). 
From the implicit function theorem, we have 
\be{dxdk1}
\frac{dx}{dk}=-\frac{\frac{\partial F}{\partial k}}{{\frac{\partial F}{\partial x}}}=
-\frac{\tan(x) -\frac{x}{\cosh^2(kx)}}{\frac{k}{\cos^2(x)}-\frac{k}{\cosh^2(kx)}}. 
\ee
 Replacing this into (\ref{derivata}), multiplying it by $c$ and expressing everything in terms of $k$ we obtain (after multiplication by $k^2$) 
 \be{derivative1}
 \frac{k^2}{1+k^2} \frac{ds}{dc}=\frac{\tan(x) -\frac{x}{\cosh^2(kx)}}{\frac{k}{\cos^2(x)}-\frac{k}{\cosh^2(kx)}} (1+k^2)-k^2x. 
\ee
 From standard trigonometric and hyperbolic identities we obtain that the denominator that appears in (\ref{derivative1}) is 
\be{stan}
 \frac{1}{\cos^2(x)} -\frac{1}{\cosh^2(kx)} =\tan^2(x)+\tanh^(kx)>0. 
\ee
Therefore the sign of $\frac{ds}{dc}$ coincides with the sign of 
\be{deriv2}
(1+k^2) \left( \tan(x)- \frac{x}{\cosh^2(kx)} \right)- k^2 x \left( \tan^2(x)+ \tanh^2(kx) \right). 
\ee
Now, in (\ref{deriv2}),   replace $\frac{1}{\cosh^2(kx)}$ with $1-\tanh^2(kx)=1-k^2 \tan^2(x)$ because of (\ref{Ffan}) and  $\tan^2(x)+\tanh^2(kx)$ with $(1+k^2)\tan^2(x)$. 
The expression simplifies to 
$
(1+k^2)(\tan(x)-x), 
$
which,  again because of (\ref{Ffan}), has the same sign as $\tanh(kx)-kx$    which is negative according to the FACT above recalled. 

\vspace{0.25cm} 

Therefore the function $s=s(c)$ is monotone and it has a limit as $c \rightarrow 0^{+}$ and $c \rightarrow 1^{-}$. In particular since $s(c)> \frac{\pi}{c}$, we have 
$\lim_{c \rightarrow 0^+} s(c)=\infty$. Furthermore, denoting by $\tilde s$ the $\lim_{c \rightarrow 1^-} s(c)$ we have that $\pi \leq \tilde  s  \leq \frac{3\pi}{2}$ and from  (\ref{sintc}) $\tilde  s$ satisfies 
$$
\tan(\tilde  s)=\lim_{c \rightarrow 1^-} \frac{1}{\sqrt{1-c^2}}\tanh(\sqrt{1-c^2} \tilde  s)=\tilde  s. 
$$

\subsubsection*{ Function $x_{int}=x_{int}(c)$}

We want to show that $|x_{int}(c)|$ is e decreasing function of $c$. Equivalently 
we can consider $|x_{int}|$ as a function of $k$ and show, using (\ref{dkdc}),  
that it is an {\it increasing} function of $k$. Therefore we consider the function 
$\frac{\sinh(\sqrt{1-c^2} s)}{\sqrt{1-c^2}}$ written in terms of $x=cs=x(k)$ and $k$, that is,  the function, 
\be{gk}
g(k)=\frac{\sqrt{1+k^2}}{k} \sinh(kx(k)), 
\ee
and prove that $\frac{dg}{dk} > 0$. Using the chain rule, we have, 
\be{prifor}
\frac{dg}{dk} =-\frac{1}{k^2 \sqrt{1+k^2}} \sinh(kx)+ \frac{\sqrt{1+k^2}}{k} \cosh(kx)\left( x+ k \frac{dx}{dk} \right). 
\ee
The derivative $\frac{dx}{dk}$ was already computed in (\ref{dxdk1}) and using (\ref{stan}) and (\ref{Ffan}) and replacing into (\ref{prifor}), we get 
$$
\frac{dg}{dk} =-\frac{1}{k^2 \sqrt{1+k^2}} \sinh(kx)+ \frac{\sqrt{1+k^2}}{k} \cosh(kx)\left( 
x-\frac{\left( \tan(x)-\frac{x}{\cosh^2(kx)} \right)}{(1+k^2)\tan^2(x)} \right), 
$$
which, using $\frac{1}{\cosh^2(kx)}=1-\tanh^2(kx)$ and (\ref{Ffan}), becomes after simplifications, 
$$
\frac{dg}{dk} =\frac{1}{k^2 \sqrt{1+k^2} \cosh(kx)} \left( 
- \tanh(kx) + \frac{kx(1+ \tan^2(x))- k \tan(x)}{\tan^2(x)}\right). 
$$
The quantity in parenthesis, using (\ref{Ffan}), and standard trigonometric identities is 
$$
-k \tan(x)+\frac{kx}{\sin^2(x)} -\frac{k}{\tan(x)}= \frac{kx}{\sin^2(x)} -k \left( \frac{\tan^2(x)+1}{\tan(x)} \right)= \frac{1}{\sin^2(x)} \left( kx -k \tan(x)\right)>0, 
$$ 
using (\ref{Ffan}) and the FACT recalled above. This completes the proof that $\frac{dx_{int}}{dk}>0$.

\vspace{0.25cm}

\subsection*{Case $1<c<\frac{2}{\sqrt{3}}$}

\subsubsection*{Function $s_{int}=s_{int}(c)$}

As before,  we define $x=cs$, so that from $\frac{\pi}{c} < s < \frac{2\pi}{c}$, we have $\pi< x< 2 \pi$. Furthermore, we 
define $k:=\frac{\sqrt{c^2-1}}{c}$, with $0< k < \frac{1}{2}$,  and $s_{int}$ which  for simplicity we denote by  $s$,  is defined implicitly by ($y=0$ in (\ref{yt}) (\ref{cass2})) (cf. (\ref{Ffan}))
\be{implidef}
F(k,x)=k \cos(kx) \sin(x)-\sin(kx)\cos(x)=0. 
\ee
Similarly to (\ref{derivata}), we have, from $x=cs$, and using the chain rule and $\frac{dk}{dc}= \frac{1}{c^2 \sqrt{c^2-1}}$, 
$$
\frac{ds}{dc}=\frac{1}{c} \left(  \frac{dx}{dc} -s \right) = \frac{1}{c^2} \left( c \frac{dx}{dk} \frac{dk}{dc} -x \right)= 
\frac{1}{c^2} \left( \frac{dx}{dk} \frac{1}{c \sqrt{c^2-1}} -x \right). 
$$
Writing everything in terms of $k$, we obtain 
\be{intermsof}
\frac{ds}{dc}= (1-k^2) \left(  \frac{(1-k^2)}{k}\frac{dx}{dk} -x \right). 
\ee
The computation of $\frac{dx}{dk}$ is done using the implicit function theorem as in (\ref{dxdk1}). We calculate 
\be{parFparx}
\frac{\partial F}{\partial x} = (1-k^2) \sin(kx) \sin(x), 
\ee
and 
\be{parFpark}
\frac{\partial F}{\partial k}= \cos(kx) \sin(x)- kx \sin(kx) \sin(x) - x \cos(kx) \cos(x).  
\ee
Replacing  $\frac{dx}{dk}$ in (\ref{intermsof}), we obtain, after simplifications, 
$$
\frac{ds}{dc}=\frac{1-k^2}{k^2 \sin(kx) \sin(x)} \left( kx \cos(kx) \cos(x) -k \cos(kx) \sin(x) \right) , 
$$
which,  using (\ref{implidef}),  becomes 
\be{intermsof2}
\frac{ds}{dc}=\frac{1-k^2}{k^2 \sin(kx) \sin(x)} \left( kx \cos(kx)-\sin(kx)\right) \cos(x). 
\ee
Now since $\pi < x < 2\pi$ and $0 < k < \frac{1}{2}$, and therefore $0 < kx < \pi$, the first factor in (\ref{intermsof2}) is always negative and the second term is always negative as well. Therefore, the sign of $\frac{ds}{dc}$ is negative (positive) when $\cos(x)$ is negative (positive), that is,  when $\pi < x < \frac{3\pi }{2}$ ($\frac{3\pi}{2}  < x < 2\pi$). The minimum is obtained when $x=cs=\frac{3\pi}{2}$. Replacing this in (\ref{implidef}), we obtain that the minimum is obtained when $k=\frac{1}{3}$ which gives $c=\frac{3}{2 \sqrt{2}}$ and $s=\pi \sqrt{2}$. 

Similarly to what was done for the case $|c|<1$, the function $s_{int}$ is monotone (and bounded) 
for $c \rightarrow 1^+$ and therefore it has a limit. If we denote the limit by $\tilde s$, from (\ref{yt}),  (\ref{cass2}),  it has to satisfy 
$$
\tan (\tilde s)=\lim_{c \rightarrow 1^+} \frac{\tan(\sqrt{c^2 -1} \tilde s)}{\sqrt{c^2-1}}=\tilde s. 
$$ 
At the point $c=\frac{3}{2 \sqrt{2}}$ which is $k=\frac{1}{3}$, the function is continuous according to the implicit function theorem and the derivative (\ref{parFparx}).

\subsubsection*{Function $x_{int}=x_{int}(c)$}
To show that the function $|x_{int}|=|x_{int}|(c)$ is decreasing as a function of $c$, we  show that it  is decreasing as a function of $k:=\frac{\sqrt{c^2 -1}}{c}$ since $k$ is increasing as a function of $c$. To do that,  we show that the function $\frac{1}{c^2 -1} \sin^2(\sqrt{c^2-1} s_{int}(c))$ is decreasing as a function of $k$. Write such a 
function in terms of $k$ and $x=cs$, 
as $\frac{1-k^2}{k^2} \sin^2(kx(k))$, and since $0 < kx<\pi$, the sin is always positive and we can consider the function (cf. (\ref{gk})),   
$$
g(k)=\frac{\sqrt{1-k^2}}{k} \sin(k x(k)). 
$$
We have 
\be{tobeput}
\frac{dg}{dk}=\left[ \frac{d}{dk} \left( \frac{\sqrt{1-k^2}}{k}  \right)\right] \sin(k x(k))+ 
\frac{\sqrt{1-k^2}}{k}  \frac{d}{dk} \sin(k x(k))=-\frac{1}{k^2 \sqrt{1-k^2}} \sin(kx)+
\frac{\sqrt{1-k^2}}{k}  \frac{d}{dk} \sin(k x(k)). 
\ee
The calculation of $ \frac{d}{dk} \sin(k x(k))$ is done using the chain rule and the implicit function theorem (first equality in (\ref{dxdk1}))  with (\ref{parFparx}) (\ref{parFpark}). We obtain 
$$
 \frac{d}{dk} \sin(k x(k))= \frac{\cos(kx(k))}{(1-k^2) \sin(kx) \sin(x)} \left[ x \sin(kx) \sin(x)-k \cos(kx) \sin(x)+ kx \cos(kx) \cos(x)  \right]. 
$$
Replacing this in (\ref{tobeput}), we obtain 
$$
\frac{dg}{dk} =\frac{1}{k^2 \sqrt{1-k^2} \sin(kx) \sin(x)} \left[ -\sin^2(kx) \sin(x) + kx \cos(kx) \sin(kx) \sin(x) -k^2 \cos^2(kx) \sin(x)+ k^2x\cos^2(kx) \cos(x) \right]. 
$$
Using (\ref{implidef}), this becomes 
\be{firstp}
\frac{dg}{dk} =\frac{1}{k^2 \sqrt{1-k^2} \sin(kx) \sin(x)}  \left[ \sin^2(kx) + k^2 \cos^2(kx) \right] \left[ x \cos(x) -\sin(x) \right]. 
\ee
Since $\pi < x< 2\pi$ and $0 <  kx < \pi$  the first factor in (\ref{firstp}) is negative while the second factor is clearly positive. The third factor can be written, using (\ref{implidef}), as 
$$
\left[ kx \cos(kx) - \sin(kx) \right] \frac{\sin(x)}{\sin(kx)}, 
$$
which is positive since $w \cos(w)-\sin(w)$ is negative for $0 < w < \pi$. 
Therefore $|x_{int}|$ is a strictly decreasing function of $k$ and therefore of $c$. Using (\ref{rc2f}) and the limit for $s_{int}$ as $c\rightarrow 1^+$, one obtains the limit 
$\lim_{c \rightarrow 1^+} |x_{int}(c)|= \sqrt{1 + \bar s}$, and using the limit   
$c \rightarrow \frac{2}{\sqrt{3}}^-$, of 
$s_{int}$, one obtains that this limit for $|x_{int}|$ is $1$.

\section*{Appendix C: Proof of Proposition \ref{intersect}}

In the main text, we have described four types of geodesics according to the value of $c$. The ones for $0< |c|< 1$, say Type 1, the ones for $|c|=1$, say Type 2, the ones for $1 < |c|\leq  \frac{2}{\sqrt{3}}$, say Type 3, the ones for 
$|c| > \frac{2}{\sqrt{3}}$, say Type 4. Geodesics of type 3 and 4 are described by the same formula but are distinguished by the fact that the ones of type 3 end on the $x$ axis while the ones of type 4 end on the unit circle. 
For all types of geodesics let us denote by $x_{int}=x_{int}(c)$ the $x$ coordinate of the terminal point, whether it belongs to the $x$-axis or the unit circle (or both). The discussion in the main text and in Appendix B has demonstrated that $x_{int}$ is a strictly increasing function of $c$, for $0< c < \infty$. Therefore among the geodesics for $c>0$ there exists only one from the point $(1,0)$ to a given point on the $x$ axis 
(with $-\infty <x \leq  -1$) or on the unit circle (with $-1 \leq x  < 1$). Recall that for  points on the unit circle (including the initial point $(1,0)$) the metric is not defined and they do not belong to the manifold  we are considering which is the plane with the unit disk removed. 

Now assume that two geodesics $\gamma_{c_1}$ and $\gamma_{c_2}$ with $0 < c_1 < c_2$ (and therefore $x_{int}(c_1) > x_{int} (c_2)$) intersect. The intersection has to be transversal, that is, the tangent vectors at the intersection point cannot coincide because this would contradict the uniqueness of the solution of differential equations since both geodesics satisfy the geodesic (second order) equations.  
The computation (\ref{dotx}) (\ref{doty}) shows that the slope of the geodesic $\gamma_c$ (near $s=0$)  is $\frac{\dot y}{\dot x}(s)=\tan(cs)$. Furthermore the slope varies as $\frac{d}{ds}  \frac{\dot y}{\dot x}(s)= \frac{d}{ds}  \tan(cs)=\frac{c}{\cos^2(cs)}$, which gives $c$ at $s=0$. Therefore the geodesic corresponding to $c_1$ is more `flat' than the one corresponding to $c_2$ and ends at a value of $x_{int}$ which is less than the value for $c_2$. If there is a point of (transversal) intersection, there 
must be  at least two because otherwise we would have   $x_{int}(c_1)    \geq  x_{int}(c_2)$.

Consider then two such points $P_1$ and $P_2$ with $P_1$ occurring before $P_2$. We know 
that $\gamma_{c_2}$ is a length minimizing geodesics between $P_1$ and $x_{int}(c_2)$ (it is, in fact, length minimizing between each point $\gamma_{c_2}(t_0)$ and $x_{int}(c_2)$, for $ 0 < t_0 < s_{int}(c_2)$) and $\gamma_{c_1}$ is length 
minimizing  between $P_1$ and $x_{int}(c_1)$.  However, this implies that both $\gamma_{c_1}$ and 
$\gamma_{c_2}$ are length minimizing between $P_1$ and $P_2$. This implies that the curve that goes from $P_1$ to $P_2$ following, say, $\gamma_{c_1}$ and then switches to, say, $\gamma_{c_2}$ is length minimizing from $P_1$  to, say, $x_{int}(c_2)$. The only points where two geodesics can and do intersect is on the $x$ axis where the two geodesics corresponding to $c$ and $-c$ for $0< |c|  \leq  \frac{2}{\sqrt{3}}$ meet. 

\section*{Appendix D: Lie algebra automorphisms of $sl(2)$ (Proof of Lemma \ref{autom})}

Let $\alpha$ be a Lie algebra  automorphism of $sl(2)$ and $[\alpha]$ the corresponding matrix in the basis $\{A_0, A_1, A_2\}$. The automorphism condition $\alpha([A,B])=[\alpha(A), \alpha(B)]$ for every $A,B \in sl(2)$ can be written in terms of the adjoint representation as 
\be{relas1}
[\alpha] ad_A=ad_{\alpha (A)} [\alpha]. 
\ee
From (\ref{commurel}), we have
\be{addav}
ad_{A_0}=\frac{1}{2} \begin{pmatrix} 0 & 0 & 0 \cr 0 & 0 & 1 \cr 0 & -1 & 0 \end{pmatrix}, \qquad ad_{A_1}=\frac{1}{2}\begin{pmatrix} 0 & 0 & 1 \cr 0 & 0 & 0 \cr 1 & 0 & 0 \end{pmatrix}, \qquad ad_{A_2}= \frac{1}{2}\begin{pmatrix} 0 & -1 & 0 \cr -1 & 0 & 0 \cr 0 & 0 & 0 \end{pmatrix},  
\ee
and writing 
$$
[\alpha]:=\begin{pmatrix} a_{00} & a_{01} & a_{02} \cr 
a_{10} & a_{11} & a_{12} \cr 
a_{20} & a_{21} & a_{22} \end{pmatrix}, 
$$
we can write (\ref{relas1}) as, for $j=0,1,2$ 
$$
[\alpha] ad_{A_j}=\sum_{k=0}^2 a_{jk} ad_{A_k} [\alpha]. 
$$
From this, using the expressions (\ref{addav}),  we obtain constraints on the entries $a_{jk}$ of $[\alpha]$ which are given by 
$$
[\alpha]=\begin{pmatrix} C_{00} & -C_{01} & -C_{02} \cr 
- C_{10} & C_{11} & C_{12} \cr 
-C_{20} & C_{21} & C_{22} \end{pmatrix}, 
$$
where $C_{jk}$ is the {\it cofactor}  corresponding to $a_{jk}$. Let $C$ denote the matrix of cofactors. 
The above relations give 
\be{toberep}
[\alpha]=I_{1,2} C I_{1,2}. 
\ee 
This implies $\det(C)=\det\left( [\alpha]  \right)$. Furthermore, using the fact that  $[\alpha]^{-1} =\frac{1}{\det\left([\alpha] \right)} C^T$, we obtain $\det\left[\alpha]\right)=1$. Therefore $[\alpha]^{-1}=C^T$ and we have $[\alpha]=I_{1,2} [\alpha]^{-T} I_{1,2}$. Then we have 
$$
[\alpha]^TI_{1,2} [\alpha]=\left( I_{1,2} [\alpha]^{-1} I_{1,2} \right) I_{1,2} [\alpha]= I_{1,2}. 
$$
Therefore $[\alpha] \in SO(1,2)$. Viceversa, if $[\alpha]  \in SO_0(1,2)$ from $[\alpha]^TI_{1,2} [\alpha]=I_{1,2}$, we get $I_{1,2} [\alpha]^T I_{1,2}=[\alpha]^{-1}=C^T$ and therefore (\ref{toberep}). 
\end{document}